\documentstyle[fullpage,12pt,amssymb]{amsart}
\def\SBIMSMark#1#2#3{
 \font\SBF=cmss10 at 10 true pt
 \font\SBI=cmssi10 at 10 true pt
 \setbox0=\hbox{\SBF Stony Brook IMS Preprint \##1}
 \setbox2=\hbox to \wd0{\hfil \SBI #2}
 \setbox4=\hbox to \wd0{\hfil \SBI #3}
 \setbox6=\hbox to \wd0{\hss
             \vbox{\hsize=\wd0 \parskip=0pt \baselineskip=10 true pt
                   \copy0 \break%
                   \copy2 \break% 
                   \copy4 \break}}
 \dimen0=\ht6   \advance\dimen0 by \vsize \advance\dimen0 by 8 true pt
                \advance\dimen0 by -\pagetotal
 \dimen2=\hsize \advance\dimen2 by .25 true in
  \openin2=publishd.tex
  \ifeof2\setbox0=\hbox to 0pt{}
  \else 
     \setbox0=\hbox to 3.1 true in{
                \vbox to \ht6{\hsize=3 true in \parskip=0pt  \noindent  
                \input publishd.tex 
                \vfill}}
  \fi
  \closein2
  \ht0=0pt \dp0=0pt
 \ht6=0pt \dp6=0pt
 \setbox8=\vbox to \dimen0{\vfill \hbox to \dimen2{\copy0 \hss \copy6}}
 \ht8=0pt \dp8=0pt \wd8=0pt
 \copy8
 \message{*** Stony Brook IMS Preprint #1, #2 ***}
}

\newtheorem{thm}{Theorem}[section]
\newtheorem{cor}[thm]{Corollary}
\newtheorem{lem}[thm]{Lemma}

\theoremstyle{remark}

\numberwithin{equation}{section}
\newcommand{\thmref}[1]{Theorem~\ref{#1}}

\newcommand{\lemref}[1]{Lemma~\ref{#1}}
\newcommand{\corref}[1]{Corollary~\ref{#1}}

\theoremstyle{definition}

\newcommand{\QED}{\rlap{$\sqcup$}$\sqcap$\smallskip}

\newcommand{\di}{\partial}
\newcommand{\dibar}{\bar\partial}
\newcommand{\ra}{\rightarrow}

\def\ssk{\smallskip}

\def\lqq{\lq\lq}
\def\sm{\smallsetminus}
\def\bolshe{\succ}
\def\ssm{\smallsetminus}
\def\tr{{\text{tr}}}

\newcommand{\diam}{\operatorname{diam}}
\newcommand{\dist}{\operatorname{dist}}

\newcommand{\inter}{\operatorname{int}}
\renewcommand{\mod}{\operatorname{mod}}
\newcommand{\tl}{\tilde}

\newcommand{\orb}{\operatorname{orb}}

\newcommand{\id}{\operatorname{id}}

\newcommand{\meas}{\operatorname{meas}}

\newcommand{\Dil}{\operatorname{Dil}}
\newcommand{\Ker}{\operatorname{Ker}}
\newcommand{\tg}{\operatorname{tg}}
\newcommand{\codim}{\operatorname{codim}}
\newcommand{\isom}{\approx}
\newcommand{\esssup}{\operatorname{ess-sup}}

\newcommand{\SLa}{\underset{\La}{\Subset}}

\def\loc{{\mathrm{loc}}}

\newcommand{\eps}{\epsilon}
\newcommand{\la}{{\lambda}}
\newcommand{\La}{{\Lambda}}
\newcommand{\Om}{{\Omega}}
\newcommand{\om}{{\omega}}

\newcommand{\AAA}{{\cal A}}
\newcommand{\BB}{{\cal B}}
\newcommand{\CC}{{\cal C}}
\newcommand{\DD}{{\cal D}}
\newcommand{\EE}{{\cal E}}
\newcommand{\EEE}{{\cal O}}
\newcommand{\II}{{\cal I}}
\newcommand{\FF}{{\cal F}}
\newcommand{\GG}{{\cal G}}
\newcommand{\JJ}{{\cal J}}
\newcommand{\HH}{{\cal H}}
\newcommand{\KK}{{\cal K}}
\newcommand{\LL}{{\cal L}}
\newcommand{\MM}{{\cal M}}
\newcommand{\NN}{{\cal N}}

\newcommand{\QQ}{{\cal Q}}
\newcommand{\RR}{{\cal R}}
\newcommand{\SS}{{\cal S}}
\newcommand{\TT}{{\cal T}}
\newcommand{\UU}{{\cal U}}
\newcommand{\VV}{{\cal V}}
\newcommand{\WW}{{\cal W}}
\newcommand{\XX}{{\cal X}}
\newcommand{\YY}{{\cal Y}}
\newcommand{\ZZ}{{\cal Z}}

\newcommand{\C}{{\Bbb C}}
\newcommand{\D}{{\Bbb D}}
\newcommand{\Hyp}{{\Bbb H}}

\newcommand{\N}{{\Bbb N}}

\newcommand{\R}{{\Bbb R}}
\newcommand{\T}{{\Bbb T}}
\newcommand{\V}{{\Bbb V}}
\newcommand{\U}{{\Bbb U}}
\newcommand{\W}{{\Bbb W}}

\newcommand{\Z}{{\Bbb Z}}

\newcommand{\tT}{{\mathrm{T}}}
\newcommand{\tD}{{\mathrm{D}}}

\def\Bf{{\bold{f}}}
\def\Bg{{\bold{g}}}

\def\Bh{{\bold{h}}}

\def\Bj{{\bold{j}}}

\def\B0{{\bold{0}}}

\newcommand{\QL}{{\cal{QL}}}

\renewcommand{\lq}{``}
\renewcommand{\rq}{''}

\catcode`\@=12

\def\Empty{}
\newcommand\oplabel[1]{
  \def\OpArg{#1} \ifx \OpArg\Empty {} \else
  	\label{#1}
  \fi}
		
\newcommand{\comm}[1]{}
\newcommand{\comment}[1]{}

\begin{document}

\bigskip\bigskip

\title[Regular or Stochastic]{Almost every real quadratic map
\\ is either regular or stochastic}
\author {Mikhail Lyubich }
%\thanks{
%   This work was supported in part by Sloan Research Fellowship
% and NSF grants DMS-8920768 and DMS-9022140.}
%\date{July 11, 1997}

\maketitle
\thispagestyle{empty}
\SBIMSMark{1997/8}{July, 1997}{}

\begin{abstract} We prove 
%the Renormalization Conjecture on 
uniform hyperbolicity of the  renormalization operator for all possible real
 combinatorial types. We derive from it that the set of 
   infinitely renormalizable parameter values in the real quadratic family
 $P_c: x\mapsto x^2+c$ has zero measure. This yields the statement
  in the title (where \lq{ regular}" means to have an attracting cycle  and \lq{stochastic}"
  means to have an absolutely continuous invariant measure). 
An application to the MLC problem is given. 
% We also improve our previous result on real parameter values where MLC holds. 
\end{abstract}

%\tableofcontents

\section{Introduction}

The main goal of this paper is to prove 
hyperbolicity of the renormalization operator  for all possible real combinatorial types, 
and to derive from it the well-known Regular or Stochastic Conjecture. 

The Renormalization Conjecture was stated by Feigenbaum \cite{F1,F2}
 and independently by Coullet \& Tresser \cite{CT,TC}
in 1978, for the particular case of doubling combinatorics. In the works of
Lanford, Epstein, Eckmann, Sinai,  Sullivan, McMullen  \cite{La1,E,EE,VSK,S1,S2,McM2},
among many others, a spectacular progress in this problem has been achieved
(see \cite{universe} for more historical details and references).
 However, until recently only the doubling case had been essentially resolved. 
 In \cite{universe} the Conjecture was proven for bounded combinatorics.
% hyperbolicity of the renormalization has been proven in
In this work we will extend the conjecture 
(compare with  Lanford's Conjecture \cite{La2}  for circle maps)
 and  prove  it for all possible combinatorial types.

 Let us introduce a few notations and state the results. 
Let us consider the real quadratic family $P_c: x\mapsto x^2+c$, $c\in [-2, 1/4]$
(in this parameter  range the maps $P_c$ have an invariant  interval).
A map $f=P_c$ is called renormalizable if it has a periodic interval $L\ni 0$ of some 
\lq{renormalization period}"
$p>1$, i.e., $f^p L\subset L$. The set of renormalizable quadratics is the union of 
disjoint closed intervals $J_k\subset (-2, 1/4)$ called \lq{renormalization windows}"
(see e.g., \cite{MvS}). The renormalization period $p_k=p(J_k)$ is constant through the window.

Let $\QL_\R$ stand for the space of real quadratic-like maps, and let $\HH_\R(f)\subset \QL_\R$ 
stand for the real  hybrid class via $f\in \QL$ (see \S \ref{q-l maps} for the definitions).
By a renormalization strip $\TT_{J_k}\subset \QL_\R$ we mean the union   of hybrid classes 
$\HH_\R(P_c)$, $c\in J_k$. On each renormalization strip one can define a real analytic 
renormalization operator $R_{J_k}: \TT_{J_k}\ra \QL_\R$, so that $R f$ is an
 appropriately restricted $p_k$-fold  iterate of $f$.
These maps can be organized in a single piecewise analytic map
$R: \cup \TT_{J_k}\ra \QL_\R$. 

Let $\Sigma$ stand for the space of two-sided sequences of natural numbers, and
$\omega$ stand the shift on this symbolic space. We are now ready to state the Renormalization
Theorem for all real combinatorial types:

\begin{thm}[Full renormalization horseshoe]\label{renormalization}
There is a set $\AAA\subset \cup\TT_{J_k}$ so that:
\begin{itemize}
\item $\AAA$ is $R$-invariant and $R|\AAA$ is topologically conjugate to the two-sided shift $\omega$;
\item The restriction $R|\AAA$ is uniformly hyperbolic; 
\item Any stable leaf $W^s(f)$, $f\in \AAA$,  coincides with the hybrid class
  $\HH_\R(f)$ and has codimension 1; 
\item Any unstable leaf $W^u(f)$ is an analytic curve which 
  transversally passes through all real hybrid classes
    except the cusp one (corresponding to $c=1/4$);
\item The renormalization operator has
uniformly bounded non-linearity on the unstable leaves  outside a neighborhood of the cusp class;
\item The expansion factor of the branches $R_{J_k}$ goes to $\infty$ as $k\to \infty$. 
\end{itemize} 
\end{thm}

{\it Remark.}  The above contraction and expansion properties hold with respect to a
\lq{Montel metric}" induced by an appropriate Banach norm (see \S \ref{q-l maps}).
The hyperbolicity is uniform in the following sense. 
The rate of contraction on the stable foliation is uniform on a 
subset of quadratic-like maps with a definite modulus
($\mod(f)\geq \mu>0$);
the rate of expansion is uniform provided $Rf$ stays outside a neighborhood of the cusp class
(see Theorems \ref{horseshoe thm} and \ref{unstable foliation}). \medskip

A quadratic map $P_c: x\mapsto x^2+c$ 
 is called {\it regular} if it has an attracting cycle (i.e., a cycle whose multiplier has
an absolute value less than 1).  In this case, the attracting cycle is unique and attracts
almost all orbits (\cite{Singer,G-sensitive}). It
is called {\it stochastic} if it has an absolutely continuous invariant measure. 
In this case the measure is unique, weakly Bernoulli, and almost all orbits are asymptotically
equidistributed with respect to it (\cite{Le,BL2}). 
\begin{thm}[Regular or stochastic]\label{regular or stochastic}
Almost every real quadratic polynomial\\
$P_c(z)=z^2+c$, $c\in [-2,1/4]$, is either regular or stochastic.
\end{thm}

 Regular quadratic maps are also called {\it (uniformly) hyperbolic}, as they are uniformly 
 expanding outside the basin of the attracting cycle \cite{Fa,G-sensitive}. 
On the other hand, stochastic maps
can also be called {\it  (non-uniformly) hyperbolic} in the sense of the Pesin theory
(as the invariant measure automatically has a positive characteristic exponent \cite{BL2}).   
Thus one can say that {\it almost any real quadratic map is hyperbolic.}

Previously it was known that  
stochastic maps are observable with positive probability
(\cite{J}, \cite{BC}) but nowhere dense (as follows from the Yoccoz theorem, see \cite{H}).
On the other hand, the set of
 regular maps is open (obviously) and dense 
(see \cite{puzzle} for the proof of this result and
further reference comments).  Our Regular or Stochastic Theorem
completes the measure-theoretical picture of  dynamics  in  the real quadratic family.
% and confirms  the general Palis Conjecture in higher-dimensional spaces \cite{Pa}. 

Let us remind the following topological decomposition of the parameter interval (see \cite{MvS}):
$[-2,1/4]=\RR\cup\NN\cup\II$, where $\RR$ stands for the regular parameter values,
$\NN$ stands for non-regular at most finitely renormalizable parameter values, 
and $\II$ stands for infinitely renormalizable parameter values. 
The set $\SS$ of stochastic parameter values is contained in $\NN$ (this follows from
a theorem that for $f\in \II$, almost all orbits converge to an attractor of measure 0
(see \cite{G,BL,S2})). 
 Thus Theorem
\ref{regular or stochastic} will follow from the following two results:
\begin{thm}[joint with Martens \& Nowicki \cite{parapuzzle,MN}]\label{LMN} Almost every\\
non-regular real quadratic polynomial which is at most finitely renormalizable is stochastic:
  $\meas(\NN\sm \SS)=0$. 
\end{thm} 

Namely, in our joint project,
 Martens and Nowicki gave a geometric condition for
existence of an absolutely  continuous invariant measure \cite{MN}, and the author has
shown that this condition is  satisfied almost everywhere \cite{parapuzzle}. Note that it
is known that the difference $\NN\ssm \SS$ is non-empty \cite{Jo,HK,Bruin}.
\begin{thm}\label{measure 0}
The set of infinitely renormalizable real quadratics
 has zero Lebesgue measure: $\meas(\II)=0$.
\end{thm}

This result will be derived from \thmref{renormalization}. Let us give a few more
applications of that Theorem. 
For any renormalization window $J_k$, there is a canonical map $\sigma: J_k\ra [-2,1/4)$
defined as the renormalization postcomposed with the  straightening. 
Let $\{J^n_i\}$ stand for the collection of
 domains of definition of $\sigma^n$, that is the windows for the 
$n$-fold renormalization, and let 
$J^n_i(\eps)=\sigma^{-n} [-2, 1/4-\eps]\cap J^n_i$.  
\begin{thm}\label{uniformly qs}
The maps $\sigma^n: J_i^n(\eps)\ra [-2,1/4-\eps]$ are uniformly 
quasi-symmetric (with the dilatation independent of $n$ and $i$).
% depending only on $\eps$).
\end{thm}

Given a renormalizable map $f$, let $p(f)$ stand for the period of the first renormalization.
The following result improves Theorem VIII of \cite{puzzle}: 
% (where  \lq{essentially}"  is a mathematical term):  
%
\begin{thm}\label{MLC}
There is a number $\bar p$ with the following property. If $f=P_c$ is a real quadratic map
with $p(R^{n_k} f)\geq \bar p$ for  a subsequence  $n_k\to\infty$, 
then the Mandelbrot set is locally connected at $c$. Moreover, the corresponding
little Mandelbrot sets $M_{n_k}$ shrinking to
$c$ have a bounded shape.   
\end{thm}

The last statement means that the canonical homeomorphisms  of the sets $M_{n_k}$ 
onto the whole Mandelbrot set $M_*$ admit  uniformly $K$-quasiconformal extensions
to  $(\eps\diam M_{n_k})$-neighborhoods of the $M_{n_k}$ (with absolute $K$ and $\eps$). 

\ssk Let us now dwell on the main ingredients of the proof of \thmref{renormalization}.
There are three types of combinatorics to take care of: bounded, essentially bounded
and high. For bounded combinatorics, Sullivan \cite{S2} and McMullen \cite{McM2}
constructed the renormalization
horseshoe $\AAA$ and its strong stable foliation. It was proven in \cite{universe}
that the renormalization horseshoe is hyperbolic.
The idea of the proof is that in the complex analytic set up lack of hyperbolicity
yields existence of  \lq{slowly shadowing orbits}". On the other hand, such orbits 
are ruled out by the Rigidity Theorem \cite{puzzle}. 
Note that  an important part of \cite{universe}
is supplying the space of quadratic-like germs (modulo affine conjugacy) with the complex
analytic structure and demonstrating that the Douady \& Hubbard hybrid classes 
\cite{DH:pol-like} form a  foliation of the connectedness locus 
with complex codimension 1 analytic leaves.

The unbounded combinatorics can be split into two types: \lq{essentially bounded}" and
\lq{high}". In the former case, the unboundedness is produced by the saddle-node behavior of the
critical point (see \cite{puzzle,LY}).  This
phenomenon can be analyzed by means of the parabolic bifurcation theory (see \cite{D2}).
Motivated by works of A.~Epstein \cite{Ep} and McMullen \cite{McM2}, 
Ben Hinkle has proven a rigidity theorem for \lq{parabolic towers}" \cite{Hi},
geometric limits of dynamical systems generated by infinitely renormalizable maps with
essentially bounded combinatorics. 
Using this result, we prove hyperbolicity of the
renormalization operator with essentially bounded combinatorics.  Note that McMullen's argument
for exponential contraction does not seem to work in this case, and instead we make use of the
Schwarz Lemma in Banach spaces.

To treat high combinatorics we need an extensive analytic preparation 
on the geometry of the puzzle and parapuzzle which was done
in \cite{puzzle,parapuzzle}. The main geometric results of these works is linear growth
of the conformal
moduli of the \lq{principal nest}" of dynamical and parameter annuli. These
imply that the image of a renormalization 
horizontal strip of high type is a narrow \lq{vertical}" strip close to the quadratic
family. This yields strong hyperbolicity of the high type renormalization, with big
contraction and expansion factors. 
Note that it is crucial for our argument that  the results of \cite{puzzle,parapuzzle}
are proven for complex parameter values (even though in this work we are ultimately interested in
the real quadratics). 

Finally, the argument of \cite{universe} (slowly shadowing orbit versus rigidity)
glues the above ingredients together 
and yields Theorem \ref{renormalization}.

\ssk This work completes a program of study the real quadratic family
  by complex methods carried in the series of papers
   \cite{LM,attractors,puzzle,LY,parapuzzle,MN,Hi,universe}.

\ssk Let us finish with  a couple remarks concerning more general settings. 
In the one-dimensional theory there are two natural  ways to proceed:
to higher degree polynomials and to $C^2$-smooth maps.  We expect the analogous 
\lq{regular or stochastic}" statement to be  valid in generic one parameter families. 
There is still a lot  of interesting work to be done in this direction. 

One can also  formulate an analogous  conjecture for the complex quadratic family 
$z\mapsto z^2+c$. Here absolute continuity of an invariant measure can be understood
with respect to Sullivan's conformal measure on the Julia set. \lq{Almost all}" in the parameter
plane can be understood in the sense of Hausdorff dimension 
as \lq{outside a set of strictly smaller  dimension}". Of course, such a conjecture cannot
be proven prior to the MLC conjecture (though can be disproven).

A general program in real higher dimensional situation was formulated by Palis
(e.g., at the  Paris/Orsay Symposium (1995)).  Roughly speaking, it asserts that in a generic 
one parameter family, there is  typically only  finitely many attractors each of which carries
an SBR  measure and such that almost any  orbit is equidistributed with respect to one of them. 
This program initiated by the work of Benedicks \& Carleson \cite{BC2} is now being intensively
carried on (see Viana \cite{Vi}, Young \cite{Y} and further references therein).

 \smallskip{\bf Notations.} 
 $\C$, $\R$, $\Z$ and $\N$ denote as usual the complex plane, the real line, and the sets of 
 integer and natural numbers respectively;\\ 
 $\D(a,r)=\{z: |z-a|<r\}$ is the open disk of radius $r$,\\
$\D_r\equiv \D(0,r)$,  $\D\equiv D_1$; \\
$\T_r=\di \D_r$ is the circle of radius $r$, $\T\equiv \T_1$;\\
$U\Subset V$ means that $U$ is {\it compactly contained } in $V$,
 that is,  the closure $\bar U$ is compact and is contained in $V$. \\
 The closure of a set $X$ will be denoted by $\bar X$.\\
When considering the space $\C^2$, $\pi_1$ and $\pi_2$ will stand for the coordinate
projections. 
Notation $\alpha\asymp\beta$ means as usual that the ratio $\alpha/\beta$ is bounded
  away from 0 and $\infty$.\\
The quasi-conformality property will be often abbreviated as \lq{qc}". 
Similarly, \lq{qs}" will stand for \lq{quasi-symmetric}". \\
Let $$
        \Dil(h)=\esssup {\di h + \dibar h\over \di h- \dibar h}
$$
stand for the dilatation of a qc map $h$.\\
Let $P_c(z)=z^2+c$.    

\smallskip{\bf Acknowledgement.} This work was done in fall 1996 during author's visit to IHES.
It was partially written during a brief visit to  ETH in Z\"urich (January 1997). I thank both
Institutes for their hospitality. 
The main results were first announced at the analysis seminar at KTH in Stockholm (November 1996), and I wish
to thank Lennart Carleson and Michael Benedicks for  their inspiring interest. 
I also wish to thank the Stony Brook University for the 
generous sabbatical support. The work
was partially supported by the NSF grant DMS-9505833.

\section{Quadratic-like germs, puzzle and towers}\label{preliminaries}

\subsection{Space of quadratic-like germs}\label{q-l maps}
This section summarizes \cite[\S\S 3,4]{universe}. 
Recall that  a map  
$f: U'\ra U$ is called {\it quadratic-like} if it a double branched covering  between topological
disks $U,U'$ such that either $U=U'=\C$ (and then $f$ is a quadratic polynomial), 
or  $U'\Subset U$. It has a single critical
point which is assumed to be located at the origin $0$, unless otherwise is stated.

The filled Julia set  is defined as the set of non-escaping point: 
$K(f)=\{z: f^n z\in U,\, n=0,1\dots\}$. Its boundary is called the {\it Julia set},
 $J(f)=\di K(f)$. 
 The sets $K(f)$ and $J(f)$ are connected if and only if the critical point itself is non-escaping:
$0\in K(f)$. Otherwise these sets are Cantor.

The {\it fundamental annulus} of a quadratic-like map $f: U'\ra U$  is the annulus between the
domain and the range of $f$, $A=U\ssm U'$. 

 Any quadratic-like map has two fixed points counted with multiplicity. In the case of connected
Julia set these two points can be dynamically distinguished.  One of them, usually denoted by 
$\alpha$, is either non-repelling or {\it dividing}, i.e., removing of it makes the Julia set
disconnected. Another one, denoted by $\beta$, is  always {\it non-dividing}. 

For the details of the further discussion we refer to \cite{universe}.
First of all, we allow to change the domains $(U, U')$ 
of a quadratic-like map without changing \lq{its germ}" near the Julia set.
More precisely,   
let us say that two quadratic-like maps $f: U'\ra U$ and $\tl f \tl U'\ra \tl U$
represent the same {\it marked germ} if there is a  string of quadratic-like maps
 $f_k : U_k'\ra U_k$,  representing  the same germ at 0,  with  both points
$0$ and $f (0)$  contained in the same connected component $W_k$
 of $U_k\cap U_{k+1}$, $k=0,1,\dots, N-1$, and such that $f_0=f$, $f_N=\tl f$.
By \cite[\S 5.4]{McM1}, a marked quadratic-like germ has a well-defined Julia set.
% Note also that a quadratic-like germ with connected Julia set  has a unique marking.

We will consider quadratic-like germs up to {\it affine} conjugacy (rescaling),
so that near the origin they can be normalized as $f(z)=c+z^2+\dots$. Marked
quadratic-like germs modulo affine conjugacy will still be called briefly
\lq{quadratic-like germs}". We will not make notational difference between quadratic-like
germs and quadratic-like maps representing them. Note also that any quadratic polynomial
$P_c: z\mapsto z^2+c$ determines a quadratic-like germ by restricting it to a
sufficiently big round disk $\D_r$. These germs will still be called quadratic polynomials.   

Let $\QL$ stand for the space of quadratic-like germs, and $\CC$ be  its {\it connectedness
locus}, that is,  the subset of germs with connected Julia set. 
 We   supply  $\QL$ with 
topology and complex analytic structure in the following fashion.
Let $\V$ be the ordered set of topological discs $V\ni 0$ with piecewise smooth boundary, with
$U\bolshe V$ if $U\Subset V$.
 Let $\BB_V$ denote the space of normalized  analytic functions $f(z)=c+z^2+\dots$ on
$V\in\V$ continuous up to the boundary supplied with sup-norm
$\|\cdot\|_V$, and let $\BB_V(g,\eps)$ stand for the $\eps$-ball  in this space centered at $g$.

If $g\in \BB_V$ is quadratic-like on $V$ then all nearby maps $f\in \BB_V$ are quadratic-like
on a slightly smaller domain. Thus we have an embedding $\BB_V(g, \eps)\ra \QL$.  
This family of embeddings induces a
 topology and complex structure on $\QL$ (see the Appendix). 

Given a set $\XX\subset \QL$, the intersections $\XX_V=\XX\cap \QL_V$ are called the
{\it Banach slices} of $\XX$.

By Lemma \ref{properties1}, compactness in $\QL$ is equivalent to sequential compactness.
Moreover, any compact set $\KK\subset\QL$ locally sits in finitely many Banach slices
$\QL_V$ and possesses a {\it Montel metric} $\dist_M$ well-defined up to quasi-isometry.  

% All basic notions of global
% analysis on this space (analytic maps, submanifolds etc.),  are expressed in terms of its {\it
% Banach slices} $\BB_V$. 

Let $\QQ=\{P_c: z\mapsto z^2+c\}\isom \C$ stand for the quadratic family. It is a complex
one-dimensional submanifold of $\QL$. By definition, 
the Mandelbrot set $M_*\subset\QQ$ is equal to  $\QQ\cap \CC$.

Given a  marked germ $f$, let $\mod(f)=\sup \mod(A)$ where $A$ runs over the fundamental 
annuli of quadratic-like representatives of $f$.
For $\mu>0$, let $\QL(\mu,\rho)$ stand for the set of normalized quadratic-like germs with
 $\mod(f)\geq \mu$ and $|f(0)|\leq \rho$. Furthermore,
let 
$$
        \QL(\mu)=\{f\in \QL: \mod(f)\geq \mu\}.
$$
 Given a set $\XX\subset \QL$,  let $\XX(\mu)=\XX\cap \QL(\mu)$.  

\begin{lem}[Compactness]\label{compactness lemma}
For any $\mu>0$ and $\rho>0$, 
the sets $\QL(\mu,\rho)$ and $\CC(\mu)$ are compact.
Moreover, if $f_n\in \QL(\mu_n,\rho)$ with  $\mu_n\to\infty$ then the limit points of
the $f_n$ are quadratic polynomials.
\end{lem}    

\begin{pf}
      See \cite[Theorem 5.6]{McM1} and \cite[Lemma 4.1]{universe}. 
\end{pf} 

Two quadratic-like germs $f$ and $g$  are called {\it hybrid equivalent} if they are
quasi-conformally conjugate by a map $h$ with $\bar\di h=0$ a.e. on $K(f)$. 
By the Douady-Hubbard Straightening Theorem \cite{DH:pol-like}, every hybrid class $\HH(f)$
with connected Julia set intersects the quadratic family $\QQ$ at a single point $c=\chi(f)$
of the Mandelbrot set $M_*$. Thus the hybrid classes can be also labeled as $\HH_c$, $c\in M_*$.

The hybrid classes can be supplied with the Teichm\"uller-Sullivan metric (see \cite{S1}):
$$\dist_T(f,h)=\inf_h\log \operatorname{Dil}(h),$$
where $h$ runs over all hybrid conjugacies between $f$ and $g$. 
 Let us also define $\dist_{T,V}$ as a similar infimum as $h$ runs over hybrid
equivalences defined in $V$ (warning:
unlike $\dist_T$, $\dist_{T,V}$ is not a metric). 

\begin{lem}\label{Teich metric}
Let $f\in \QL_V$, $g\in \HH(f)$, and $W\Subset V$. 
There exists an $\eps>0$ such that 
if $\dist_{T,V}(f,g)<\eps$ then $g$ belongs to $\QL_{W}$ and $\|f-g\|_W<\eps$. 

Vice versa, for any $\eps>0$ there is a $\delta>0$ such that if $g\in\BB_V(f,\delta)$
then $\dist_{T,W}(f,g)<\eps$. 
\end{lem} 

\begin{pf} The first statement follows from the fact that any normalized qc map $h$ 
with $\Dil(h)<\eps$ is uniformly close to $\id$.

To prove the second one, observe that if  $g\in\BB_V(f,\delta)$, then  $f$ and $g$ have
$(1+\eps)$-qc equivalent fundamental annuli $A_f$ and $A_g$ in a slightly smaller domain, such that
the corresponding qc map $h: A_f\ra A_g$ respects dynamics on the inner boundaries of the annuli.
Such an $h$ extends to a hybrid equivalence between $f$ and $g$ in $W$ with dilatation $(1+\eps)$.
\end{pf} 

It is proven in \cite{universe} that the hybrid classes $\HH_c$, $c\in M_*$, are
connected complex codimension one submanifolds of $\QL$.  They form a foliation $\FF$
(or rather a \lq{lamination}")  called {\it horizontal}. This foliation is
transversally quasi-conformal everywhere, and complex analytic on $\inter\CC$.  

Let us state the former result more precisely. Take two hybrid equivalent maps
$f_i$, and two complex analytic transversals $\SS_i$ to the leaf $\HH\equiv \HH(f_i)$
 via $f_i$. The holonomy $\gamma$ from $\SS_1$ to $\SS_2$ along $\FF$ 
is called {quasi-conformal} if  
 it admits a qc extension to a neighborhood of the $f_i$ in the transversals
 (beyond the connectedness locus).  

The local dilatation of $\gamma$ at $f_1$
is 
$$
   \inf_h \lim_{\eps\to 0} \Dil (h\, |\, \D(f_1, \eps)),
$$
where the infimum is taken over all local  qc extensions $h$ of $\gamma$.

 \begin{thm}[\cite{universe}, Theorem 4.14]\label{trans qc}
 Given two quadratic-like maps as above, the holonomy $\gamma$ is quasi-conformal. 
 If $\mod(f_i)\geq\mu$ then the local  dilatation of $\gamma$ at $f_1$ is bounded by $K(\mu)$.
 Moreover, if additionally $\dist_T(f_1,f_2)\leq \rho<1$ then 
 the above dilatation is $1+O(\rho)$ with the constant depending on $\mu$ only.
 \end{thm}
 
\comm{
 \begin{pf} The maps $f_1$ and $f_2$ are included in a complex one-dimensional analytic
 curve $\XX\subset \HH$, the {\it Beltrami path} $f_\la$. This path belongs to finitely
 many Banach slices $\BB_{V_i}$. Within each of these slices the holonomy $\gamma$ is
 quasi-conformal by the $\la$-lemma (see the Appendix). Hence it is quasi-conformal altogether.
 
 % Moreover, by compactness of $\QL(\mu)$, there is a $\bar\rho=\bar\rho(\mu)$ 
 % such that if $\rho<\bar\rho$ then the maps $f_i$ belong to the same Banach slice $\BB_V$
 Moreover, the hyperbolic distance between the $f_i$ in $\XX$ is $O(\rho)$.
 By the $\la$-lemma, the local dilatation of $\gamma$ is $O(\rho)$ as well.
 \end{pf} 
end comm}

Let $\EE$ denote the space of real analytic expanding circle endomorphism  of degree
two considered up to rotation. By definition, any $g\in \EE$ admits a
double covering  complex analytic
extension $g: V\ra V'$ to symmetric annuli neighborhoods of the circle $\T$
(with piecewise smooth boundary)
 such that $V\Subset V'$ (where \lq{symmetric}" is understood in the sense of the
involution $\tau: \bar\C\ra \bar\C$ about the circle $\T$).
%  Note that such a map can be normalized so that $g(1)=1$.

There is a projection 
\begin{equation}\label{pi}
           \pi: \QL \ra \EE,
\end{equation}
 which associates to an $f\in \QL$
its external map $g=\pi(f)\in \EE$ (see   \cite{DH:pol-like}, \cite[\S 3.2]{universe}).
The construction goes as follows. Take a quadratic-like representative  
$f: V'\ra V$ and consider the fundamental annulus $A=V\sm V'$, $\mod(A)=\mu$. 
 Using the map $f: I\ra O$
from the inner to outer boundary of $A$, attach to the inner boundary of $A$
 an abstract annulus $A_1$ of modulus $\mu/2$. It comes together with a double covering
$A_1\ra A$ which extends $f: I\ra O$. Using this covering, attach  in a similar way an annulus 
$A_2$ of modulus $\mu/4$ to the inner boundary of $A_1$, etc. 
Taking the infinite union of these annuli together with
$\C\sm V$, we obtain  a conformal punctured disk $S$ and a double covering $F$
between annuli neighborhoods of its ideal boundary. Let us
uniformize it, $\phi=\phi_f: S\ra\C\sm \D$, and conjugate $F$ by $\phi$. This provides us with
a double covering between  outer annuli neighborhoods  of $\T$. Reflecting it about the circle,
we obtain the desired external map $g: V'\ra V$,
 where $V$ and $V'$ are  symmetric annuli neighborhoods of $\T$ and $V'\Subset V$. 

As the external map $g$ is defined up to rotation so that it can be normalized by
putting its fixed point at 1: $g(1)=1$. 

Note further that the uniformization $\phi_f$ provides a conformal isomorphism between  the
fundamental annuli $U\sm U'$ and $V\sm  (V'\cup \D)$, and conjugates $f$ and $g$
on their inner boundaries. Moreover, by means of dynamics it can be analytically extended
to a domain  containing the
critical value $f(0)$ (and containing the critical point 0 on its boundary). 
Thus we can consider the image of the critical point under this map:
\begin{equation}\label{xi}
    \xi(f)=\phi_f (f(0))
\end{equation}
(it is well-defined once $g$ is normalized). 
The inverse map $\psi_f=\phi_f^{-1}$ will be called the {\it uniformization of $f$ at}
$\infty$.

Let us consider a real Banach space $\BB^s_V$ of $\T$-symmetric (i.e., commuting with $\tau$) 
complex analytic maps $V\ra \C$  continuous in $\bar V$. 
A sufficiently small Banach neighborhood $\BB^s(f,\eps)$ consists of expanding circle 
endomorphisms, so that it is embedded into $\EE$.
This endows $\EE$ with the inductive limit topology and real analytic structure.

\comm{ Moreover, in  \cite{universe} $\EE$ is supplied with a \lq{Bers-Sullivan}" 
{\it complex} analytic structure.  
 There is a projection
$\pi: \QL\ra\EE$   which associates to any $f\in \QL$ its {\it external map}. 
end comm}
Restricted to any hybrid class $\HH_c$, $c\in M_*$, the projection $\pi$ becomes a homeomorphism.
The inverse map $i_c: \EE\ra \HH_c$ is provided by the \lq{mating}" of a circle map
$g\in \EE$ with the quadratic polynomial $P_c$ 
(see   \cite{DH:pol-like,universe}). 
 This allows us to transfer the complex analytic structure  from the hybrid class $\HH_0$ of
$z\mapsto z^2$ to the space $\EE$. 

The Bers-Sullivan complex structure makes
the projection $\pi:\QL\ra\EE$ and all the parametrizations $i_c: \EE\ra \HH_c$, $c\in M_*$,
 complex analytic (see \cite[\S 4.3]{universe}). The fibers $\ZZ_g$, $g\in \EE$,
 of $\pi$ turn out
to be complex analytic curves in $\QL$ \cite[Theorem 4.18]{universe}.
 They are called {\it vertical fibers}. 

\comm{
 The external map $g=\pi(f)$ admits a natural normalization so that $g(1)=1$. 
The maps $f$ and $g$ are conformally equivalent in appropriate annuli neighborhoods:
there exists a conformal map $\phi_f: U'\sm \Omega\ra V'\sm \La$ conjugating $f| U\sm \Omega$
to $g| V\sm \La$. Here $\Omega=\Omega_f$ and $\La=\La_f$ 
are topological  disks containing respectively
$K(f)$ and $K(g)=\D$, and moreover $\Omega$ contains the critical value $c=f(0)$.
This provides us with a special local chart for $f$ which will be called its
{\it external model} (see \cite[\S 3.3]{universe}).

This allows us to extend the Riemann mapping  $\C\sm M_*\ra \C\sm \bar\D$ 
of the complement of the Mandelbrot set to a smooth  mapping
\begin{equation}\label{xi}
   \xi: \QL\sm\CC\ra\C\sm \D
\end{equation}
to the complement  of the connectedness locus.  It is defined as the position of the
critical value $f(0)$ in the external model of $f$, $\xi(f)=\phi_f(f (0))$. 
Moreover, $\xi$ is {\it vertically holomorphic}, i.e., it holomorphic on the vertical fibers
(\cite[Lemma 4.9]{universe}).
end comm}

The map $\xi$ (\ref{xi}) provides a smooth  extension (actually, real analytic)
 of the Riemann mapping  
$\C\sm M_*\ra \C\sm \bar\D$ (see \cite{Orsay})
to the complement  of the connectedness locus. 
Moreover, this map is {\it vertically holomorphic}, 
i.e., it holomorphic on the vertical fibers $\ZZ_g$
\cite[Lemma 4.9]{universe},

Note finally that the Green function $G=\log |\xi|: \QL\sm \CC\ra \R_+$ 
provides us with a dynamically natural  way to measure the \lq{distance}"  from an
$f\in \QL\ssm \CC$ to the connectedness locus.
The level sets of the Green function are called
{\it equipotentials}  $\EEE_r$ (of radius $r>1$)  in $\QL$.  
% as the pull-back of the
% round circle $\{z: |z|=r\}$ under the  map $\xi$.
 One can show that they are real
codimension one smooth submanifolds in $\QL\sm\CC$ (since $\xi$ is a smooth submersion).

\comm{ The map $\QL\ssm\CC\ra \EE\times(\C\ssm M_*) $, $f\mapsto (\pi(f),\xi(f))$,
is a local  smooth diffeomorphism, so that the  both components
$\pi$ and $\xi$ are submersions.  
Let us define the {\it equipotential}  $\EEE_r$ of level $r>1$  in $\QL$
as the pull-back of the round circle $\{z: |z|=r\}$, under the  map $\xi$. 
Since $\xi$ is a submersion,  the equipotentials  are
the codimension one smooth submanifolds in $\QL\sm\CC$. Moreover,
the \lq{Green function}" $G(f)=\log\xi(f)$ 
continuously extends to the connectedness locus as 0.  Let
$\DD_r=\{f: G(f)< \log r\}$ be the domains bounded by the equipotentials of level $r$. 
This gives us a dynamically natural  way to measure the \lq{distance}"  from an
$f\in \QL\ssm \CC$ to the connectedness locus. 
end comm}

\subsection{Quadratic-like families}\label{q-l families} 
% We assume that the reader is familiar with the theory of holomorphic motions
% (see \cite{survey,MSS,Sl}).
The reader is referred to \cite{DH:pol-like,parapuzzle} 
for a  discussion of quadratic-like families.
Let us have a domain $\La\Subset \C$.
 A domain $\V\subset \La\times\C$ is called a {\it tube} over 
$\La$ if it is homeomorphic over $\La$ to a straight tube $\La\times \D$.
Let  $V_\la=\pi_1^{-1} \{\la\}$ stand for the vertical fibers of a tube $\V$.
We will  assume that  they are bounded by piecewise smooth curves and contain 0.  
Let $\B0=\La\times \{0\} $ stand for the zero section of $\V$. 
If $\La'\subset \La$ then  let $\V|\La'=\V\cap (\La'\times\C)$.
Let $\di^h \V=\cup_{\la\in \La} \di V_\la$ stand for the {\it horizontal boundary} of $\V$. 

Given two tubes $\V'\subset \V$ as above, let us say that $\V'$ is compactly contained in $\V$
over $\La$, $\V\underset{\La}{\Subset} \V'$, if   
the relative closure of $\V'$ in $\La\times \C$
   is bounded,   contained in $\V$,  and $\dist(\di^h \V', \di^h \V) >0$.

By definition, a map $\Bf: \V'\ra \V$ between two tubes $\V'\subset \V$ over $\La$
     is called a {\it quadratic-like family} over  $\Lambda=\La_\Bf$ if
 $\Bf$ is  a holomorphic  endomorphism preserving the fibers, 
  and such that every fiber restriction $f_\la: V_\la'\ra V_\la$, $z\mapsto z^2+c(\la)+\dots$, 
   is a normalized quadratic-like map with a critical point at 0. 
Clearly  any quadratic-like family represents a holomorphic curve in $\QL$.

Let 
$$
      M_\Bf=\{\la: 0\in K(f_\la)\}=\Bf\cap \CC
$$ 
stand for the Mandelbrot set in $\Bf$. The straightening provides  a canonical continuous map 
\begin{equation}\label{chi}
           \chi: M_\Bf\ra M_*
\end{equation}
 (see Douady \& Hubbard \cite{DH:pol-like}).
  The family $\Bf$ over $\La$ is called {\it full} if $M_\Bf\Subset \La$.
 In this case the straightening properly maps $M_\Bf$ onto the whole Mandelbrot set $M_*$.
A full family is called {\it unfolded} if  the straightening (\ref{chi}) is injective 
(and hence is a homeomorphism).

Let  $\phi(\la)=f_\la(0)$ stand for the critical value of $f_\la$, and let
$\Phi(\la)=(\la,\phi(\la))$ stand for the \lq{critical value section}" $\Lambda\ra \V$. 
The family $\Bf$ is called {\it proper} if $\V'\SLa\V$ and the map $\Phi$ is proper over $\La$, 
i.e, for any $\W\SLa \V$, we have $\Phi^{-1} \W\Subset \La$. 
% $\phi(\la)\in \di V_\la$ for $\la\in\di \Lambda$. 
In particular, for a proper family, the domain
 $$
        \Lambda'\equiv \La'_\Bf=\Phi^{-1}\V'=\{\la: f_\la(0)\in V_\la'\}
$$
is compactly contained in $\La$, so that this family is full.  

For a full family, one
defines the {\it winding number} $w(\Bf)$ as the winding number of the curve 
$\la\mapsto\phi(\la)$ about the origin, 
as $\la$ goes once anti-clockwise around a Jordan curve $\Gamma$ enclosing $M_\Bf$. 
The family is unfolded if and only if  it has winding number 1.
In this case there is a single  superstable parameter value  $*$
 (the root of $\phi$) called the {\it center} of  $\Lambda$.  

%   By the Douady and
% Hubbard theorem \cite{DH:pol-like}, if $\Bf$ is proper and unfolded then $M_\Bf$ is
%canonically homeomorphic to the standard Mandelbrot set $M_*$.   
% Let  $\La''=\{\la\in \La': f_\la^2(0)\in V_\la'\}.$

Given a proper unfolded quadratic-like family $\Bf=\{f_\la\}$, let 
$$
  \mod(\Bf) =\min[\;\mod(\Lambda\sm \Lambda'),\; \mod(\La'\sm M_\Bf),\; 
         \inf_{\la\in \Lambda}\mod(V_\la\sm V_\la')\;],
$$ 
where $\mod(\La\sm \La')$ is  understood
%(in case these domains are not simply connected)
as the extremal length of the family of curves in $\La\sm \La'$
joining $\di \La$ with $\di \la'$, and similarly for $\mod(\La'\sm M_\Bf)$.

Let $\GG$ stand for the collection of proper unfolded
 quadratic-like families up to affine change of variable in $\la$.
Such a family can be normalized so that superattracting parameter value $*$ sits at the origin
and $\diam M_\Bf=1$.  
% conformal change of variable in
% $\la$ (so that we can always assume that $(\Lambda,*)=(\D_r,0)$). 
We will impose the {\it Carath\'eodory} topology on $\GG$ (compare \cite{McM1}). 
In this topology a sequence of families
$\Bf_n$ over $\Lambda_n$ converges to a family $\Bf$ over $\Lambda$ if
$(\Lambda_n, *)$ Carath\'eodory converges to $(\Lambda,*)$ 
%(after appropriate reparametrization),
the  domains $(\V_n, \V_n',\B0)$ (with the preferred zero sections) Carath\'eodory converge to
$(\V,\V',\B0)$, so that the convergence respects the fiberwise structure, and finally the  $\Bf_n$
converge to $\Bf$ uniformly on compact sets. 

Carath\'eodory convergence in the above definition means:
\begin{itemize}
\item For any tube $\W\ni \B0$ compactly contained in $\V_n$,
 all the domains $\V_n$ eventually contain $\bar \W$;
\item Vice versa: Any domain $\W$ as above compactly  contained in infinitely
  many of the $\V_n$, is also contained in $\V$.
\end{itemize}
Note that 
the Carath\'eodory convergence of the families corresponds to the uniform on compact sets
convergence of the corresponding holomorphic  curves in $\QL$.

Let $$
  \GG_\mu=\{\Bf\in \GG:  \diam \V\leq \mu^{-1},\; \mod(\Bf)\geq\mu\},
$$
where the family $\Bf$ above is meant to be normalized:  $\diam M_\Bf=1$. 

\begin{lem}\label{compactness}
 For any $\mu>0$, the space $\GG_\mu$ is compact.
\end{lem}

\begin{pf} Let us have a sequence of normalized families $\Bf_n$ over $(\La_n,0)$. 
Since $\mod(\La\sm \La')\geq \mu>0$, $\mod(\La'\sm M_\Bf)\geq \mu$, and
$\diam M_\Bf=1$,  the families of domains $\La_n\equiv \La_{\Bf_n}$ and
$\La'_n\equiv \La'_{\Bf_n}$ are Carath\'eodory compact. Select a
Carath\'eodory converging subsequence: 
$(\Lambda_n,0)\to (\Lambda,0)$ and  $(\Lambda_n',0)\to (\Lambda',0)$.
 Since the $\mod(\La_n\ssm \Lambda_n')$ stay away from 0 and the $\diam \La_n$ are bounded, 
the limit domain $\Lambda'$ is compactly contained in $\La$. 

Take a domain $\Om\Subset \La$. 
Let us now select a converging subsequence of the domains $(\V_n,\B0)\to (\V,\B0)$ over $\Om$. 
To this end let us consider a family  $\WW$ of domains $\W\ni {\bold 0}$ 
such that some relative neighborhood of $\bar \W$ in $\Om\times \C$ is
compactly contained in infinitely many of the $\V_n$.
This family is non-empty. Indeed,  by normalization of the $f_\la$ and the bound
$\mod(f_\la)\geq\mu>0$, $\WW$ contains the round tubes $\Om\times \D_\eps$
with sufficiently small $\eps>0$.  

The family $\WW$ has a countable basis $\WW^0$, 
i.e. a countable family of domains such that for any $\W\in \WW$ there exists a $\W^0\in \WW^0$
such that $\W^0\subset \W$ 
(for instance, take polygonal domains in $\WW$ with rational vertices).  
Select an exhausting sequence $\W_n\in \WW^0$ so that $\W_1\subset \W_2\subset\dots$ and
no $\W\in \WW^0$ contains all the $\W_n$. Now select a
sequence $\V_n$ such that any $\W_n$ is compactly contained in all the $\V_k$, $k\geq n$. 
This subsequence Carath\'eodory converges to $\W=\cup\W_n$ over $\Om$. 

Now select an
exhausting sequence of domains $\Om\Subset\La$. By means of the diagonal procedure 
we can find a  subsequence of tubes $\V_n$ converging over each $\Om$.
Hence by definition it is Carath\'eodory  converging over the whole $\La$,
$\V_n\to \V$. 

Similarly  select a further subsequence so that  $\V_n'\to \V'$. 

Furthermore, since the maps $\Bf_n=\{f_\la\}$
 act fiberwise as branched double coverings and normalized as
$z\mapsto z^2+c(\la)+\dots$ at the origin,
 by the Koebe Theorem they form a pre-compact family on each compact subset of
$\V_n'$. Hence we can select a converging subsequence $\Bf_n\to \Bf$.
As the limit family is clearly proper and unfolded,  we are done. 
\end{pf}

Let $f\in\CC$ with $\chi(f)=c$,
 and let $P_f=\tD(i_c\circ \pi)_f$ be the projection of the tangent space
$\tT_f\QL$ onto the tangent space $\tT_f\HH_f$ to the leaf.  
  The one-dimensional spaces 
\begin{equation}\label{kernels}
K_f=\operatorname{Ker} \tD\pi_f=\Ker\tD P_f
\end{equation}
 form a continuous subbundle $\KK\subset \tT\QL$ 
complementary to the tangent subbundle $\tT\FF$. For a vector $u\in \tT_f\QL$, let
$u^h=P_f u$ and $u^v=(I-P_f)u$ stand for its \lq{horizontal}" and \lq{vertical}"
projections. 

 If we take a 
 Banach slice $\QL_V$ whose tangent space contains $K_f$, then for $u\in \tT_f\QL_V$ we can
define the angle between $u$ and $\tT_f\HH_f$ by letting 
\begin{equation}\label{theta}
\tg(\theta)=\|u^v\|/\|u^h\|.
\end{equation}
 Small angle means that $v$ is almost tangent to  the leaf.  

If we have a family of tangent vectors belonging to finitely many Banach slices, we say
that they are {\it uniformly transversal} to the foliation $\FF$ if their angles with the
foliation stay away from 0.

\begin{lem}\label{transversality}
Any family $\Bf\in \GG_\mu$ is uniformly transversal to the foliation $\FF$ with the
lower bound on the  angle depending only on $\mu$.
\end{lem}

\begin{pf} 
Otherwise, by compactness (Lemma \ref{compactness}),
 there would be a family $\SS$ in $\GG_\mu$ which were
 tangent to some leaf  $\LL$ of $\FF$. Let us take a Banach slice  $\QL_V$ locally containing 
$\SS$ and a curve $\gamma\subset \LL$ tangent to $\SS$. Then the slice $\LL_V$ is still
tangent to $\SS$ in $\QL_V$. 

By \cite[Lemma 4.12]{universe},
the Banach slice $\FF_V$ of the foliation
$\FF$ is a Banach foliation with codimension 1 complex analytic leaves . 
 Moreover, it is transversally analytic over the $\inter \CC_V$. 
Let us   apply to it the results of the Appendix, \S \ref{function theory}. 

 By the Hurwitz Theorem, 
$\SS$ would have the same number of intersection points (counted with multiplicity) with all 
nearby leaves of $\FF_V$. But by the Intersection Lemma, the intersection points with the
nearby leaves of $\inter \CC_V$ are simple, so that there would be more than one such a point.
 But unfolded families intersect every leaf of $\FF$ at a 
single point (\cite{DH:pol-like}) - contradiction.
\end{pf} 

Let us now show that the quadratic-like families $\Bf\in \GG_\mu$ \lq{uniformly overflow}"
the connectedness locus. This can be measured in terms of the function 
$\xi: \QL\ssm\CC\ra \C\ssm\bar\D$ (\ref{xi}), and means that any family in question \lq{goes
beyond}" an equipotential $\EEE_r$ with $r=r(\mu)$. 

\begin{lem}\label{uniform overflowing}
There is an $r=r(\mu)>1$ such that for any family $\Bf$ over $\Lambda$,
$|\xi(f_\la)|>r$ for $\la\in\di \Lambda$. 
\end{lem}

\begin{pf}  Let $f_\la: U_\la'\ra U_\la$ and $g_\la: V_\la'\ra V_\la$ be the corresponding
external map. Since  and $f_\la 0\in U_\la\sm U_\la'$ for $\la$ near $\di \La$ 
and $\xi(f_\la)$ corresponds to $f_\la(0)$ in the external model, we have:
$\xi(f_\la)\in V_\la\sm V_\la'$. 
Since 
$$
    \mod(V_\la\ssm V_\la')=\mod(U_\la\sm U_\la')\geq \mu,
$$ 
there is an annulus of modulus at least $\mu/2$ which separates  $\xi(f_\la)$
  from the unit circle. 
Hence $|\xi(f_\la)|>r(\mu)>1$ for $\la$ near $\di \La$, and we are done.
\end{pf} 

If $\mod(\Bf)$ is big then the family $\Bf$ is close to the quadratic family in the following 
sense:

\begin{lem}\label{close to quadratic}
For any $\eps>0$ and $r$, there is a $\mu$ and a Banach space $\BB_V$ containing disk
$\D_r=\{c: |c|<r\}$ in the quadratic family $\QQ$ with the following property.
 If $\Bf$ is a full unfolded family with 
$\mod(\Bf)>\mu$, then there is topological disk $\Delta\subset \Bf$ which belongs to
$\BB_V$ and is represented in that space as a graph of an analytic function $\phi: \D_r\ra E$
(where $E$ is a complement of $\QQ$ in $\BB_V$) with $\|\phi\|<\eps$.
 \end{lem} 

\begin{pf} Take a domain $V$. If $\mu$ is big enough then all 
quadratic-like maps $f\in \QL(\mu)$ clearly  belong to a  Banach slice $B_V$.
Select a complement $E$ to $\QQ$ in this Banach space, and let $p: \BB_V\ra \QQ$ be the
projection of $\BB_V$ onto $\QQ$ parallel to $E$. 
By Lemma~\ref{compactness lemma},
 \begin{equation}\label{small difference}
\|P_{c(\la)}-f_\la\|_V=\eps<\eps(\mu),
\end{equation}
 where $P_{c(\la)}=p(f_\la)$ and $\eps(\mu)\to 0 $ as $\mu\to \infty$. 

 Let us take a big $\rho$ and consider the curve 
 $\gamma=\{f_{\la(t)}\}$ in $\Bf$ parametrized
by a Jordan curve $\delta\subset\La$ close to $\di\Lambda$ 
(where $\Bf$ is a family over $\Lambda$). 
Since $\Bf$ is unfolded, the winding number of $\Gamma=p\circ\gamma$ around 0 is equal to 1.

Moreover,  by (\ref{small difference}),  $\Gamma$ encloses a disk $\D_r\subset \QQ$ with
a big $r=r(\mu)$. Indeed, if $\delta$ is sufficiently close to $\di\La$, then
the critical value $f_\la(0)$ is arbitrary close to $\di U_\la$ for $\la\in \delta$.
Hence $f_\la(0)$ can be separated from the Julia set $J(f_\la)$ by a fundamental
 annulus of modulus at least $\mu/2$. Hence the critical value $P_{c(\la)}0$ is separated
from the Julia set $L(P_{c(\la)})$ by a fundamental annulus of modulus at least
$\mu/2-\delta(\eps)$ where $\delta(\eps)\to 0$ as $\eps\to 0$. It follows that
$|c(\la)|\geq r(\mu)\to \infty$ as $\mu\to \infty$.

Hence the winding number of $\Gamma$ about the disk $\D_r$ is equal to 1. 
It follows that the projection of $\Bf$ onto $\QQ$ univalently covers this disk.
Together with (\ref{small difference}) this yields  the desired statement. 
\end{pf}

Let us say that a Mandelbrot set $M_\Bf$ has a $K$-{\it bounded shape} if it is
 canonically homeomorphic to the standard set $M_*$ by a map which admits a $K$-qc extension
 to a neighborhood $D$ of $M$ with $\mod(D\ssm M)>1/K$.  Theorem \ref{trans qc} yields:
 
 \begin{lem}\label{bounded shape}
 For any $\Bf\in \GG_\mu$, the Mandelbrot set $M_\Bf$ has  $K(\mu)$-bounded shape.
 \end{lem}

Let us say that a quadratic-like family $\Bf: \V'\ra \V$ over $(\La, *)$ is {\it equipped} if
it is supplied with a holomorphic motion $\Bh$ of the fundamental annulus  over $\La$, 
 $h_\la: V_*\sm V_*'\ra V_\la\sm V_\la'$, $\la\in \La$. 
For instance, the quadratic family is equipped over a domain  $\La$ bounded by any parameter
equipotential.  
All primitive Mandelbrot copies
in a full equipped family $(\Bf, \Bh)$ are generated by equipped
quadratic-like families (see \cite{parapuzzle}). 

Given a holomorphic motion $\Bh$ over $\La$, let 
   $$\Dil (\Bh)=\sup_{\la\in \La} \Dil(h_\la).$$
For an equipped quadratic-like family $(\Bf, \Bh)$, let us call the pair of numbers
$(\mod(\Bf), \Dil(\Bh))$ its {\it geometry}. We say that the geometry is bounded by 
over a collection of equipped families, if $\mod(\Bf)\geq \mu>0$ and
$\Dil(\Bh)\leq K$ for all families of the collection. 

Let $\GG_\mu^e$ denote the collection of equipped quadratic-like families $(\Bf, \Bh)$
with $\mod(\Bf)\geq \mu$, $\Dil(\Bh)\leq \mu^{-1}$. 

% Now we are going to equip the external fibers $\ZZ_g$. 
For $f\in \QL$, let us consider  the following objects:
\begin{itemize}
\item the uniformization  $\psi_f$ defined after (\ref{xi});
\item the projection $\Pi: \QL\ra \HH_0$, $\Pi=(\pi|\HH_0)^{-1}\circ \pi$, where
   $\pi$ is the projection (\ref{pi}) and $\HH_0$ is the hybrid class of $P_0: z\mapsto z^2$;
\item for $G\in \HH_0$, we will use the notation $R_G\equiv\psi_G: \C\sm \D\ra \C\sm K(G)$  
     for the Riemann   mapping with positive derivative at $\infty$;
\item the map $\Psi_f=\psi_f\circ R_{G}^{-1}$, where $G=\Pi(f)$, which is another conformal
  representation of  $f$ \lq{near $\infty$}" .
\end{itemize}

\begin{lem}\label{uniformization}
  The conformal representation $\Psi_f(z)$ analytically depends on $f$.
\end{lem}

\begin{pf} By \cite[Lemma 4.10]{universe}, 
\begin{equation}\label{composition}
     (f,z)\mapsto \Psi_f(z)
\end{equation}
 is continuous.
In particular, this means that if  $\Psi_{f_0}$ is defined on a domain $\Omega_0$ and
$\Omega\Subset \Omega_0$, then all nearby maps $f\in\UU\equiv \BB_U(f_0,\eps)$ 
are defined on $\Omega$. Let us  show that (\ref{composition}) is analytic on 
$\UU\times\Omega_0.$

Let us consider a holomorphic family $f_\la\in \UU$ of quadratic-like maps over $(\La, 0)$.
For $\la$ near the origin we can select  a fundamental annulus $A_\la$ holomorphically moving
with $\la$ in such a way that  the corresponding holomorphic motion $h_\la: A_0\ra A_\la$ respects
the boundary dynamics (see \cite[Prop. 9]{DH:pol-like} or \cite[Lemma 4.2]{universe}). Let us
consider the corresponding holomorphic  family of conformal structures
$\mu_\la=h_\la^* (\sigma)$ on $A_0$, where $\sigma$ is the standard structure on $A_\la$.  Pulling
them back by the external map $g_0$, we obtain a holomorphic family of conformal structures on the
Riemann surface $S_0\isom \C\sm \bar\D$, ($\mu_\la$  is extended on $\C\ssm U_0$ as the standard
structure). Let us extend these structures to $\bar \D$ as the standard ones as well. We obtain a
holomorphic family of complex structures on $\C$
 which  will be still denoted as  $\mu_\la$.

By the Measurable Riemann Mapping Theorem, there is family of qc maps $\om_\la$ holomorphically
depending on $\la$ which solves the Beltrami equations $(\omega_\la)_*\mu_\la=\sigma$. It maps
$\C\sm \bar\D$ onto $\C\sm K(G)$ where $G=\Pi(f)\in \HH_0$.
% Let us consider the Riemann mapping $R: \C\sm \bar\D \ra  \C\sm K(G)$. 
Then $\Psi_\la=h_\la\circ \om_\la^{-1}$,
 and we conclude that it analytically depends on $\la$.  
\end{pf}

% We say that a functional on $\QL$ is {\it vertically holomorphic} if it is holomorphic
%  on the vertical fibers $\ZZ_g$, $g\in \EE$.
%
% \begin{cor}\label{psi}
% The uniformization $\psi_f$ depends on $f$  smoothly and vertically holomorphic.
% \end{cor}

For $G\in \HH_0$, let $\ZZ_G=\Pi^{-1} \{G\}$ stand for the corresponding vertical fiber.
Given a Banach neighborhood $\VV$ in $\HH_{0,V}$, the set
 $$\TT_\VV=\cup_{G\in \VV} \ZZ_G\subset \QL$$
will be called  a {\it vertical tube over $\VV$}. 
Let us say that a tube $\TT=\TT_\VV$ is {\it equipped} if 
there is a  neighborhood $\UU\subset\TT$ of  $\TT\cap \CC$ with the following properties
\begin{itemize}
\item The vertical fibers $\UU_G=\UU\cap \Pi^{-1} G$ represent full unfolded quadratic-like
    families;
\item There is a  holomorphically moving fundamental annulus 
$A_{G,\la},\; G\in \VV, \la\in \UU_G$  such that the motion respect the boundary dynamics.
\end{itemize}
In other words, we equip all the vertical fibers $\ZZ_G$ over $\VV$ in the way analytically
depending on $G\in \VV$.

\begin{lem}\label{equipment}
Any   $G_0\in \EE$ belongs to an equipped  vertical tube $\TT_V\ni G_0$. 
% Any vertical fiber $\ZZ_g$, $g\in \EE$,
% represents a full unfolded holomorphic family which can
%be equipped near its Mandelbrot set $M_g$. 
%Moreover, the geometry of the corresponding equipped
%family depends only on $\mod(g)$. 
\end{lem}

\begin{pf} External fibers represent full unfolded quadratic-like families 
 by \cite[Corollary 4.19]{universe}. In order to equip them, let us use the  
conformal representations $\Psi_f$. By Lemma~\ref{uniformization},
 $\Psi_f$  holomorphically depends on $f$.
% Let $G_0=(\pi|\HH_0)^{-1} g$. 
% Let us select a fundamental annulus $A_0=V_0\sm V_0'$ on $\C\sm K(G_0)$

By \cite[Lemma 4.2]{universe}, 
there is a Banach neighborhood $\VV=\HH_{0,V}(G_0, \eps)$ which can be equipped 
with a holomorphically moving fundamental annulus
$A_G=V_G\sm V'_G\subset \C\sm K(G)$. Moreover, the conformal representation 
$\Psi_f$ is well defined on $\C\sm V_G'$, $G=\Pi(f)$, provided $f$ varies in a certain
neighborhood $\UU\subset \TT_\VV$ of the connectedness locus $\CC_\VV=\TT_\VV\cap \CC$. 
It follows that $B_f\equiv U_f\sm U_f'=\Psi_f(A_G)$, $G=\Pi(f)$,  is a
holomorphically moving fundamental annulus of $f\in \UU$.
% over the region $\Lambda\subset \ZZ_g$ where 
% $\psi_f| A$ is well-defined. So $\La$ is equipped  with this holomorphic motion $\Bh$. 

 Furthermore, $\UU$ is the union of $\CC_\VV$ and  a domain  where
 $a(f)\equiv \Psi_f^{-1}(f(0))\in V$.
By \cite[Theorem 3.4]{universe}, for any $a\in V\sm K(G)$ and $G\in \VV$,
 there exists a unique $f=\theta(G,a)\in \ZZ_G$ such that $\Pi(f)=G$ and
 $a(f)=a$ (\lq{mating}" of $G$ and $a$). 
% Moreover, by \cite[Lemma 4.9]{universe}, this map is conformal. 

It follows that  $\UU_G\sm \CC$ is a
topological annulus, so that $\UU_G$ is a topological disk representing an unfolded
quadratic-like family. 
If $f\in \di \UU_G$ then by definition $a(f)\in \di V$ and hence $f(0)\in \di U_f$.
Thus this family is proper.  
\comm{
To control the geometry of this equipped family we may need to shrink the domain $\La$ a bit. 
Consider, e.g., the domain $\La'=\{f\in \ZZ_g: \xi(f)\in V'\}$. Since $\xi$ conformally maps
$\La\sm \La'$ onto $V\sm V'$, $\mod(\La\sm \La')=\mod(V\sm V')\equiv \mu.$
By the $\la$-lemma, $\Bh$ has a bounded dilatation depending  only on $\mu$. 
Moreover, by construction, $\mod(U_f\sm U_f')=\mu$ for all $f\in \La$. 
Hence the geometry of $(\Bf, \Bh)$ is controlled by $\mu$ only. Selecting the fundamental
annulus $V\sm V'$ to be almost optimal (i.e, with almost maximal possible modulus),
we obtain the last statement.   }  
\end{pf}

\subsection{Puzzle, parapuzzle and renormalization}\label{puzzle sec}
The notion of {\it complex renormalization} was introduced by Douady and Hubbard 
\cite{DH:pol-like,D} in order to explain computer observable little Mandelbrot copies
inside the Mandelbrot set (see \cite{M,McM1,puzzle} for further discussion).
 
Let $f$ be a quadratic-like map. Assume that we can find topological disks $U'\Subset U$
around 0 and an integer $p$ such that $g= f^p: U'\ra U$ is a quadratic -like map with
connected Julia set. Assume also that the \lq{little Julia sets}" $f^k J(g)$, $k=0,\dots, p-1$,
 are pairwise disjoint except, perhaps, touching at their non-dividing $\beta$-fixed points. 
  Then
the map $f$ is called renormalizable (with period $p$) and the quadratic-like germ $g$ considered
up to rescaling is called a \lq{renormalization}"  of $f$. The map $f$ can be renormalizable
with different periods, finitely or infinitely many. Accordingly it is called 
{\it at most finitely } or {\it infinitely renormalizable}.  There is a canonical way
to produce the {\it first} renormalization $Rf$ of $f$, with the smallest period. It is provided
by the {\it Yoccoz puzzle}.  

The reader can consult \cite{puzzle}, \S 3,  for a detailed discussion of the combinatorics
of the Yoccoz puzzle. The main combinatorial object considered in that work is the
{\it principal nest} of puzzle pieces $V^0\supset V^1\supset\dots$. There is a flexibility
of the choice of the first puzzle piece $V^0$. For the sake of this work (focused on the
real combinatorics) it can be selected as the Yoccoz puzzle piece $Y^{(1)}$
of the first level (bounded by the external rays landing at the dividing fixed point  $\alpha$,
the symmetric point $\alpha'$, and some equipotential). 

Then $V^{n+1}$ is inductively defined as the pull-back of $V^n$ corresponding to the first
return of the critical point back to $V^n$. The corresponding return map $g_n: V^n\ra V^{n-1}$
is a branched double covering. The return to level $n-1$ is called {\it central} if 
$g_n(0)\in V^n$. Let $n_k$ count the non-central levels. If this sequence is infinite
then the map $f$ is non-renormalizable. 
Otherwise the principal nest ends up 
with an infinite central cascade $V^{n-1}\supset V^{n}\supset\dots$, and the map
$g_n: V^n\ra V^{n-1}$ (after perhaps little thickening of the domain and the range) is
a quadratic-like map with connected Julia set. The germ of this map (up to rescaling) 
is called the first renormalization $Rf$ of $f$.  

Let us  now state a little lemma which will be useful in what follows:

\begin{lem}\label{little} 
No quadratic polynomial $P_c$ can be realized as the renormalization $Rf$
of a quadratic-like map. 
\end{lem}

\begin{pf} Indeed, the renormalization $Rf$ admits the analytic continuation to the domain of
$f^p$ as a branched covering of degree $2^p>2$. This is certainly not compatible with the quadratic
extension to the whole complex plane. 
\end{pf}

The number of the non-central levels in the principal nest is called the
 {\it height } of~$f$. 

The map $g_n: V^n\ra V^{n-1}$ is a restriction of the full first return map
$g_n: \cup V^n_i\ra V^{n-1}$ (denoted by the same letter). Here $V^n_i\subset V^{n-1}$ are
puzzle pieces with disjoint interiors, $V^n_0\equiv V^n$, and the restrictions
$g_n: V^n_i\ra V^{n-1}$ are univalent for $i\not=0$. 

Let us now consider the quadratic family $P_c: z\mapsto z^2+c$. For any parameter value
$c_0\in M_*$ outside the main cardioid, there is a nest of {\it parapuzzle pieces} 
$\Delta^1(c_0)\supset \Delta^2(c_0)\supset\dots\ni c_0$ corresponding to the dynamical
principal nest. For parameter values $c\in \Delta^n(c_0)$, the \lq{combinatorics}"
of the first return maps to the puzzle piece $V^{n-1}$ stay the same (see \cite{parapuzzle}
for the precise definition which, however, does not matte for the following discussion).

If $P_c$ is non-renormalizable then the parapuzzle pieces $\Delta^n(c)$ shrink to $c$
(Yoccoz, see \cite{H}, and \cite{parapuzzle}). Otherwise the return maps
 $g_{n,c}=P_c^p: V^n\ra V^{n-1}$, $c\in \Delta^n$, on the renormalization level form a
quadratic-like family $\Bg$. In the primitive case 
(when the little Mandelbrot set $M_\Bg$ is not attached to the main cardioid,
which is equivalent to saying that  $n>2$),
$\Bg$ is a full unfolded family. In the satellite case, $\Bg$ is {\it almost}
full and unfolded which means that the straightening $\chi$ homeomorphically maps $M_\Bg$
onto \lq{unrooted}" Mandelbrot set $_*\ssm \{1/4\}$ (see \cite{D}).  

 The Mandelbrot set $M=M_\Bg$ encodes the
combinatorial data of the renormalization: all maps $P_c$ with $c\in M$ are
\lq{renormalizable with the same combinatorics}". The period of this renormalization
is certainly constant throughout the copy, $p=p(M)$. 
Moreover, the little copies produced in this way are {\it maximal} in the sense that they are
not contained in a any other copy except for the whole set $M_*$.

A Mandelbrot copy is called {\it real} if it is centered on the  real line. 
Let $\MM$ stand for the family of maximal real Mandelbrot copies, 
 We say that the maps $f\in \TT_{M}$ are renormalizable with {\it real} combinatorics.

Let $\TT_{M}\subset \QL$ stand for the set of quadratic-like germs which are hybrid equivalent
to the quadratic maps $P_c$ with $c\in {M}$ (that is, the union of the hybrid classes via $M$). 
We call it a {\it renormalization strip}. The maps in the strip are renormalizable with
the same combinatorics encoded by the little Mandelbrot set $M$.   Thus the renormalization
operator $R$ is canonically defined on the union of all the renormalization strips. 
The restriction $R|\TT_M$ will also be denoted as $R_M$. 
\begin{lem}[de Melo - van Strien \cite{MvS}]\label{injectivity} 
The renormalization operator $R:\cup\TT_\MM\ra\QL$ is injective.
\end{lem}  

Moreover, renormalization is non-expanding with respect the Teichm\"uller-Sullivan metric on
the hybrid classes:
\begin{equation}\label{non-expanding}
\dist_T(Rf, Rg)\leq \dist_T(f,g).
\end{equation}
This immediately follows from the fact that an appropriate restriction of a hybrid conjugacy
$h$ between $f$ and $g$ provides a hybrid conjugacy between the renormalizations 
$Rf$ and $Rg$. This observation was a starting point for Sullivan's renormalization theory
\cite{S1}. 
 
For any $M\in \MM$, there is a canonical homeomorphism $\sigma: M\ra M_*$ defined
as the composition of the renormalization and the straightening, $\sigma=\chi\circ R$
(\cite{DH:pol-like,M}). 

If a map $f\in \QL$  is renormalizable a few times, then its combinatorics is encoded by 
a sequence (finite or infinite) $\tau(f)=\{M_0, M_1,\dots\}$  such that $R^n f\in \TT_{M_n}$,
$n=0,1\dots$. One says that an infinitely renormalizable map $f$ has a {\it bounded type}
 if the periods $p(M_n)$ are bounded.  
 
Note  that any $R_M$ admits a complex analytic extension to Banach neighborhoods of
of points $f\in \TT_M$.
 Namely, if $R_M f=f^p: U'\ra U$ is a quadratic-like renormalization of $f\in \CC_V$ 
then any nearby map
$g\in \QL_V$  admits a quadratic-like return map $g^p: U_g'\ra U$ with the same range.
We can 
call this map the renormalization of $f$ even when it has Cantor Julia set.  Since $g^p$
analytically depends on $g$, this provides us with the desired extension
(see \cite[\S 5.3]{universe} for a more detailed discussion).

Let us say that a map $f$ is {\it non-escaping} under the renormalization of type
$\tau=\{M_0,M_1,\dots\}$ if all the maps  
\begin{equation}\label{f}
f_n=R_{M_n}\circ\dots \circ R_{M_0} f
\end{equation}
 are well-defined
(where the $R_{M_k}$ stand for the analytic extensions of the renormalizations)
and $\mod(f_n)\geq\eps>0$, $n=0,1,\dots$. 
\begin{lem}[\cite{universe}, Lemma 5.7]\label{escaping} 
If a quadratic-like map $f$ is non-escaping under the renormalization of type
$\tau=\{M_0,M_1,\dots\}$ then it is infinitely renormalizable with type $\tau$. 
\end{lem} 

\comm{
\begin{pf} All we need to show is that $f_n\in \CC$, where the $f_n$ are defined by (\ref{f}).

 As $\mod(f_n)\geq\eps>0$, 
there is a choice of domains and ranges 
 $f_n: U_n\ra V_n$ such  that the boundaries $\di U_n$, $\di V_n$ have
bounded geometry (i.e., they are quasi-circles with bounded dilatation).
 It follows that $\diam U_n\to 0$   (\cite{McM2}, Proposition 4.10).

Let us show that the Julia set $J(f)$ is connected. Indeed, since
$ K(f)\cap U_n\supset K(f_n)\not=\emptyset$, we have:
$\dist(0, K(f))\leq \diam U_n\to 0$. Hence $0\in K(f)$.

For the same reason, all Julia sets $J(f_n)$ are connected, and we are done.
\end{pf}
end pf}

Any equipped proper  unfolded quadratic-like family $(\Bf, \Bh)$ over a topological disk
$\La$ can be
tiled into the parapuzzles similarly to the quadratic family, depending on the initial choice of
external rays but canonical on the Mandelbrot set  (see \cite{parapuzzle}). Let $f_0\in \Bf$ be
renormalizable with type $M$. 
Then as in the quadratic case, we have a full or almost full unfolded
quadratic-like family 
$\Bg_n=\{g_{n,\la}:  V^n_\la\ra V^{n-1}_\la\}$ on the
corresponding parapuzzle $\Delta^n\ni f_0$. The renormalization 
$R_M f_\la=g_{n,\la}$  is well defined and analytic on this parapuzzle, so that this provides
us with the analytic continuation of $R_M$ onto the parapuzzle $\Delta^n$.
It will be naturally denoted as $R_M \Bf$. Moreover, this family is also equipped with a 
holomorphic motion $\Bj$, so that  we  can write under these circumstances that
  $(\Bg, \Bj)=R_M (\Bf, \Bh)$.   

Let us consider
 a conformal map $f: \SS\ra\TT$ between  two Riemann surfaces supplied
with conformal metrics. The {\it distortion} (or {\it non-linearity}) of $f$ is defined
as follows:
$$ n(f)=\sup_{z,\zeta\in \SS}\log {\|Df(z)\|\over \|Df(\zeta)\|}.$$

The following statement shows the renormalization has bounded non-linearity on unfolded
quadratic-like families:
\begin{lem}\label{nonlinearity}  Let us consider a  family
 $(\Bf: \V'\ra \V)\in \GG_\mu$ over $\Lambda$,
and let $\Lambda'=\{\la: f_\la(0)\in V_\la'\}$. Then
  $R$ has a $K(\mu)$-bounded nonlinearity on $\Bf$ over $\Lambda'$ with respect to the hyperbolic 
  metric on $\Lambda$. Moreover, the non-linearity in the hyperbolic scale $\eps$ (i.e., within
  any hyperbolic disk in $\Lambda$ of radius $\eps$) is $O(\eps)$. 
    \end{lem}

\begin{pf}%  Let us consider the set $\XX$ of proper quadratic-like families $\Bf$ with
% winding number 1 and such that
% $\mod(\Bf)\geq\nu/2$, $\Psi(f_\la)=r$. If $r$ is sufficiently small (depending only on  $\mu$)
% then $\XX$ is compact in the open compact topology.
\comm{By Lemma \ref{}, the family $\GG_\mu$ is contained in finitely many Banach slices $\BB_{V_i}$ and
is uniformly transversal to the foliation $\FF$. 
Then we can select in these slices finitely many linear coordinate charts $E^s\oplus E^n$
with $\dim E^n=1$ so that  any
$\SS\in \GG_\mu$ is a union  of graphs $E^n\ra E^s$ with bounded slopes. It follows that
the metric on $\SS$ induced from the slice is uniformly equivalent to the metric
pulled from  $E^n$.  In turn, by the Koebe Theorem, the latter metric is uniformly equivalent
to the Euclidean metric in the parameter domain $\Lambda$.}
This follows from the  Koebe Distortion Theorem.
\end{pf}

\subsection{Essentially bounded combinatorics}\label{ess bounded} 
There is a special type of renormalization combinatorics related to the parabolic 
bifurcation which usually requires a special treatment. In this section
we will describe this phenomenon. 

Let $f$ be a renormalizable map of period $p$
with real combinatorics $M\in \MM$. Let us consider a  {\it central cascade} 
\begin{equation}\label{cascade}
V^m\supset V^{m+1}\supset\dots\supset V^{m+N}
\end{equation}
meaning that all the levels $m, m+1, \dots, m+N-2$, are central:
 $g_{m+1} 0\in V^{m+N-1}\ssm V^{m+N}$. Then the quadratic-like map $g_{m+1}$ is
combinatorially close to either the Ulam-Neumann map $z\mapsto z^2-2$, or to the
parabolic map $z\mapsto z^2-1/4$ (see \cite{puzzle}, \S 8). In the former case the cascade
is called {\it Ulam-Neumann}, while in the latter it is called {\it saddle-node}. 

 For $z\in \omega(0)\cap (V^m\ssm V^{m+1})$, let us consider the level $k$ such that either
 $g_{m+1} z\in V^k\ssm V^{k+1}$ with $k<m+N$, or $g_{m+1} z\in V^{m+N}$ and then  $k=m+N$.
Let $d(z)=\max(k-m, m+N-k)$ stand for the \lq{depth}" of landing of $z$ in the cascade. 
Finally, let $d_m=\sup d(z)$ as $z$ runs through the set of points as above. 

If (\ref{cascade}) is a saddle-node cascade, then 
all the levels $m+i$ with $m+d_m<i<m+N-d_m$ will be called {\it neglectable}. 

Let us remove from the orbit $\{f^n 0\}_{n=0}^{p-1}$ all the points  whose first landing 
at some saddle-node cascade occurs on a neglectable level (we will refer to this procedure
as \lq{eliminating the neglectable part of the cascade"). The number
of points which are left is called the {\it essential period} $p_e=p_e(f)=p_e(M)$ of the
renormalization (see \cite{puzzle,LY}).  

If the period $p$ is much bigger than the essential period $p_e$ then the orbit of the
critical point spends a lot of time near a \lq{ghost}" parabolic point. Since this 
saddle-node behaviour (also called {\it intermittency}) is well-understood 
(see \cite{D2,DD}), such a combinatorial situation admits a thorough analysis.

We say that an infinitely renormalizable map has {\it essentially bounded combinatorics}
if $p_e(R^n f)\leq \bar p_e$, $n=0,1\dots$

Let us describe this phenomenon via the parameter plane. 
Let $f$ be a quadratic polynomial with  $p_e(f)\leq \bar p_e.$ 
% Let us say that a cascade (\ref{cascade}) is \lq{long}" if $N>3\bar p_e$. 
% Of course, such a cascade must  be saddle-node. 
 Let us consider the first cascade
(\ref{cascade}) with $m=0$, and the first return $f^{l_0} 0\in V^{N-1}\ssm V^{N}$ of the
critical point to this cascade. By definition of the essential period, the return time $l_0$ is
bounded it terms of $p_e$.

Take the parapuzzle piece  $\Delta^m$ corresponding to the return map $g_{1}=f^{l_0}$. This
parapuzzle contains a little Mandelbrot copy $M$ centered at the parameter value for which
$g_{1,c}0=0$. Let $b$ be the cusp of this Mandelbrot set.
If the above cascade is long enough then $f$ is a small perturbation of $f_b$ 
(moving out of the little copy $M$). 

As the return time $l_0$ is bounded, the number
 of little copies $M$ specified in this way is bounded (in terms of $p_e$). Let $\MM^0$
stand for this family of little copies. 
 
%After passing through the cascade the critical point can return a few times 
%under iterates of $g_{m+1}$ to the top levels $i\in [m, m+d_m]$ of the cascade,
%then to a bottom level $i\in [m+N-d_m, m+N-1]$. This return specifies a bounded family
%$\MM^2$ of little Mandelbrot copies near the cusps of the copies from $\MM^1$. 

Let us now wait until the first return of the critical point to the next 
cascade (\ref{cascade}) of the principal nest (with $m=n_1+1$). Let it happens at moment $l_1$.
The  combinatorics of the critical orbit until moment $l_1$ specifies finitely many 
sequences $\MM^1_i$ of little Mandelbrot copies going to the cusps of the sets from $\MM^0$. 
Namely, the length of the first cascade specifies the element of the sequence. The sequence
itself is specified by the combinatorics of the orbit after eliminating the neglectable
part of the first  cascade. 

Let us now  consider the motion of the orbit through the second cascade
until its landing at the third one. Then we obtain
finitely many Mandelbrot sequences $\MM^2_{ij}$ accumulating to the cusps of the previous
Mandelbrot sets. The sequence is specified by the combinatorics of the critical orbit 
after eliminating the neglectable parts of the first two cascades. The element of the sequence is
specified by the pair of lengths $(N_1, N_2)$ of the two cascade (so that formally speaking, it is
a double sequence). If $N_1\to \infty$ then the corresponding subsequence of the Mandelbrot copies
accumulates on the cusps of the first family $\MM^1$ (independently of $N_2$). 
If $N_1$ stays bounded while $N_2\to \infty$, then the copies accumulate on the cusps of
$\MM^2$.   

Similarly we can consider the fourth cascade and the corresponding copies accumulating on
the cusps of the first three families, etc. 

Since the height  of the principal nest is bounded in terms of  $\bar p_e$,
this procedure will terminate in a  bounded number of steps,
which altogether specifies a family  $\MM(\bar p_e)$ of little Mandelbrot copies $M$ which
contains all the copies with $p_e(M)\leq \bar p_e$
(see \cite{Hi} for a more detailed and  formal discussion of this situation).

\subsection{Geometric bounds}\label{geometric bounds}
\comm{
The following compactness lemma explains why it is important to have a lower bound on
the $\mod(f)$. Let $\QL(\eps)$ stand for the space of quadratic-like maps with
$\mod(f)\geq\eps$.  

\begin{lem}[\cite{McM1}]\label{compactness1} For any $\eps>0$, the space $\QL(eps)$  is compact.
\end{lem}
  }

An infinitely renormalizable map is said to have a priori bounds if $\mod(R^n f)\geq\eps>0$,
$n=0,1,\dots$. We say that a map $f$ is {\it close to the cusp} if it has an attracting 
fixed point with multiplier greater than $1/2$. Note that renormalizable maps are not
close to the cusp.

\begin{thm}[A priori bounds \cite{LvS,LY}]\label{a priori bounds}
Let $f$ be $n$ times renormalizable real quadratic-like map with $\mod(f)\geq \mu>0$.
Then $$\mod (R^n f)\geq \nu_n(\mu)\geq\nu(\mu)>0,$$
unless the last renormalization is of doubling type and $R^n f$ is close to the cusp. 
Moreover,  $\limsup \nu_n(\mu)\geq\nu>0$, where $\nu$ is an absolute constant.
Thus all real infinitely renormalizable maps have a priori bounds. 
\end{thm}

The following two geometric results are crucial for our study. 

\begin{thm}[Big dynamical moduli \cite{puzzle}]\label{puzzle}
 Let $\mod(f)\geq \mu>0$. Then  for any $M\in\MM$, 
$$\mod(R_M f)\geq \nu(\mu, M)\geq\bar\nu(\mu)>0,\quad f\in M, $$
unless $p(M)=2$ and $f$ is close to the cusp. 
Moreover, $\nu(\mu, M)\to\infty$ as $p_e(M)\to \infty$.
\end{thm}

{\it Remark.} A related result on moduli growth for real quadratics was independently proven
  by Graczyk \& Swiatek \cite{GS}. Note in this respect that our proof needs in a crucial way
   the above Theorem \ref{puzzle} for  complex parameter values (even
  though  in this paper  we are  ultimately interested in the real case).

\ssk The corresponding parapuzzle result is: 

\begin{thm}[Parameter moduli \cite{parapuzzle}]\label{parapuzzle}
Let $(\Bf, \Bh)\in \GG_\mu^e$ be an equipped quadratic-like family over $\La,$
% with $\mod(\Bf)\geq\mu$, $\Dil(\bh)\leq \mu^{-1}$, 
and let  $M\in \MM$, $p(M)>2$. 
 Then $(\Bg, \Bj)=R_M(\Bf, \Bh)$ is an equipped  quadratic-like family with 
$$\mod(\Bg)\geq\nu=\nu(M,\mu)\geq \nu(\mu)>0,\quad \Dil(\Bj)\leq K(\mu).$$
Moreover, for any $\mu>0$, $\nu(M,\mu)\to \infty$ as $p_e(M)\to\infty.$
\end{thm} 
{\it Remark.}
The domain of definition of the renormalized family $(\Bg, \Bj)$ can be chosen as a parapuzzle
piece $\Delta_M\subset \La$ bounded by certain parameter rays and equipotentials
\cite{parapuzzle}. Moreover,  if we equip a  vertical tube $\TT_\VV$ (see Lemma \ref{equipment}), 
then the  parapuzzle pieces  $\Delta_G$ in the vertical fibers $\ZZ_G$
will holomorphically move with $G\in \VV$. This motion is obtained by  the holonomy along
the extended foliation $\FF$. Indeed, the parapuzzle pieces are specified by the coordinate
of the critical value in the chart obtained by the straightening of the fundamental annulus.
On the other hand, by definition (see the proof of \cite[Theorem 4.13]{universe}),
 this coordinate  determines the leaf of the extended  foliation. Thus
 we obtain an
equipped tube $\cup_{G\in \VV}\Delta_G$ of the parapuzzle pieces  
to which the renormalization $R_M$ analytically extends along the vertical fibers.   
\begin{cor}\label{shrinking}
Let us consider an equipped quadratic-like family $(\Bf, \Bh)\in \GG_\mu^e$ over $\La$,
 and let $M_i=M_i(\Bf)\subset\D$ be the corresponding family of maximal real Mandelbrot copies
   except the doubling copy. 
 Then the sets $M_i$ have $K(\mu)$-bounded shape and
   $\diam(M_i)\to 0$ as $p(M_i)\to 0$ at rate depending only on $\mu$.
\end{cor} 

\begin{pf}
By  Lemma \ref{bounded shape}, the Mandelbrot set $M_\Bf$ has a $L(\mu)$-bounded shape.  
Hence it is enough to check shrinking of the $M_i$ in the case of the quadratic family $\QQ$.
Moreover by the same Lemma and  Theorem \ref{parapuzzle}, the shapes of the sets $M_i$
are bounded. Hence it is enough to have shrinking of their real traces
$M_i\cap \R$. But these traces are pairwise disjoint as the $M_i$ are maximal.
\end{pf}

{\it Remark.} The doubling renormalization ($p(M)=2$) produces an almost full unfolded
   quadratic-like family whose  Mandelbrot set misses  a single point, its cusp
   \cite{DH:pol-like}. The set left after removing a neighborhood of the cusp is
    qc equivalent to the corresponding piece of the Mandelbrot set (with dilatation
  depending only on the geometry of the family $(\Bf, \Bh)$ and the size of the
   removed neighborhood \cite{universe}).

%Let us refine this result:
%
%\begin{thm}\label{ess bounded families}
%Under the circumstances of the above lemma, there is an absolute $\nu>0$ such that
% $\nu(M,\mu)\geq\nu$ except for finitely many copies $M$ (where the exceptional family
% of copies depends on $\mu$). 
% \end{thm}

% \begin{pf} By the previous result, all the families $\Bg_i\equiv R_M\Bf$, $M\in\MM_\R$, 
% stay in a compact part of the space $\QL$, and hence the associated Mandelbrot copies
% $M_i\subset \QQ$ have a bounded shape \cite{universe}. Since they are pairwise disjoint,
% their real traces shrink. Hence $\diam M_i\to 0$ as well. 
%\end{pf}

\subsection{Renormalization limits}\label{limits}

The previous results allow us to address the problem  of possible  renormalization limits. 

\begin{thm}\label{limit points}
Let us have a sequence of real renormalizable maps
$f_k\in \TT_{M_k}$ with $\mod(f_k)\geq\mu$ and $p(M_k)\to\infty$.
 Take a limit  $f=\lim R f_k$ 
of their renormalizations. Then $f$ is either a parabolic quadratic-like map, or
a quadratic polynomial.
\end{thm} 

\begin{pf} Note first that by Theorem \ref{puzzle}, the sequence $R f_k$ is
pre-compact, so that we can always extract limit points. Furthermore, if
the essential periods $p_e(M_k)$ are uniformly bounded then the renormalizations
$R f_k$ must converge to the cusp points of the corresponding family $\MM(\bar p_e)$ of little
Mandelbrot copies (as described in \S \ref{ess bounded}). On the other hand,
if $p_e(M_k) \to \infty $ along some subsequence, then by Theorem \ref{puzzle}
the corresponding subsequence of the renormalizations must converge to a quadratic polynomial.
\end{pf}

Note that parabolic cusp points with arbitrary combinatorics can be realized as above limits.
Just take any little Mandelbrot copy $M$ and a sequence of copies $M_k$ going to its cusp $c$
produced by increasingly long parabolic cascades (keeping $p_e(M_k)$ bounded). 
Then the limit points
$R f_k$, $f_k\in M_k$,  have the combinatorics of the cusp $c$. 

Also, take any real quadratic  polynomial $P_c$ without attracting points such that $c$ is not
a parameter of doubling bifurcation. Then $P_c$ can be realized as one of the above limit.
Indeed, it is easy to construct a sequence $M_k$ converging to such a $c$,
 with $p_e(M_k)\to\infty$ (e.g., approximate $c$ by Misiurewicz points and take little 
Mandelbrot copies nearby). Then by Theorem \ref{puzzle}, $R f_k\to P_c$ for $f_k\in M_k$.

\subsection{Combinatorial rigidity} 

\begin{thm}[\cite{puzzle}]\label{rigidity}
Let $f$ and $g$ be two infinitely renormalizable  quadratic-like maps with the same real
combinatorial type $\tau=\{M_0,M_1,\dots\}$
(but not necessarily real), and 
 with {\it a priori bounds}. Then $f$ and $g$ are hybrid equivalent.
\end{thm}

Together with a priori bounds (Theorem \ref{a priori bounds}) the Rigidity Theorem yields:

\begin{cor}\label{real rigidity}
For any real combinatorial type $\tau=\{M_0, M_1,\dots\}$, there is a single
real quadratic $P_c$ with this combinatorics.
\end{cor}

\subsection{McMullen's towers}\label{McM towers}
{\it McMullen's  tower} $\bar f$
is a sequence $\{f_k\}_{k=l}^n$ of quadratic-like maps with connected Julia sets 
such that $f_k=R f_{k+1}$. Combinatorial type  $\tau(\bar f)$ 
of such a tower is a sequence  of maximal
Mandelbrot copies $M_k\in \MM$ such that $f_k\in M_k$.   Let $p(\bar f)=\sup p(f_k)$ and
$p_e(\bar f)=\sup p_e(f_k)$ stand  respectively for the \lq{period}" and the \lq{essential period}"
of the tower. One says that the tower has a bounded (or essentially bounded) combinatorics
if $p(\bar f)$ (respectively $p_e(\bar f)$) is finite. 

The modulus $\mod(\bar f)$ of the tower is defined as  $\inf\mod(f_k).$
 One says that a tower has a priori bounds if $\mod(\bar f)>0$.  The space of
towers is supplied with the weak topology: $\bar g_m=\{g_{m,k}\}_k\to \bar f$
as $m\to \infty$ if for each index $k$,
$g_{m,k}\to f_k$. Compactness of $\CC(\mu)$ yields:  

\begin{lem}\label{compactness2} The  space of towers with uniformly
bounded  combinatorics  and common a priori bounds  is compact.
\end{lem}

% The filled Julia set $K(\bar f)$ of a tower is defined as the union $\cup K(f_k)$ (it is not
% closed). Two towers are called hybrid equivalent if they are quasi-conformally conjugate by a
% map $h$ with $\dibar h=0$ a.e. of $K(\bar f)$.

%The following Towers Rigidity Theorem is very important:  

\begin{thm}[Towers rigidity]\label{towers rigidity} 
 Two  bi-infinite towers with the same bounded combinatorics and a priori bounds 
 are affinely equivalent.  
\end{thm}

\begin{pf} By the Rigidity Theorem \ref{rigidity} two bi-infinite towers with the same
combinatorics are quasi-conformally equivalent. By McMullen's Rigidity Theorem \cite{McM2}
qc equivalent towers are affinely equivalent. 
\end{pf}  

Later on we will prove a similar rigidity theorem for towers with arbitrary 
real combinatorics (see Theorem \ref{realization}).

\subsection{Parabolic towers}
Motivated by the works of C.~McMullen \cite{McM2} and A.~Epstein \cite{Ep}, 
in \cite{Hi} parabolic towers are introduced  as geometric
limits of McMullen's towers with uniformly bounded essential period. 
Fix a $\bar p_e$, and consider the family $\MM(\bar p_e)$ of little Mandelbrot copies
 associated with $\bar p_e$-essentially bounded combinatorics (see \S \ref{ess bounded}).
 A {\it parabolic tower} $\bar f$ with $\bar p_e$-essentially bounded 
combinatorics is a sequence of semigroups $\{G_n\}_{n=l}^n$ with two generators
$G_n=\{f_n, g_n\}$. The map $f_n$ is either renormalizable or parabolic. In the former case
$g_n=\id$ and $R_M f_n=f_{n-1}$ where $p_e(M)\leq \bar p_e$. In the latter case $f_n$ has
combinatorics of the cusp point of  some $M\in\MM(p_e)$, and  $g_n$
is the transit map between the Ecale-Voronin cylinders of $f_n$ (see \cite{D2,Hi}).

Two towers as  above are called combinatorially equivalent if the maps $f_n$ are either
renormalizable with the same combinatorics (i.e., labeled by the same little
Mandelbrot copy $M_n$) or parabolic with the same combinatorics. Two towers have a priori
bounds if $\mod(f_n)\geq\eps>0$, $n=0,1,\dots$.  

\begin{thm}[\cite{Hi}]\label{Ben}
If two parabolic towers $\bar f$ and $\bar g$ with $\bar p_e$-essentially bounded combinatorics    
and a priori bounds
are combinatorially equivalent then they are affinely equivalent.
\end{thm}

\section{Hyperbolicity of the renormalization operator}\label{hyperbolicity}
   
\subsection{Uniformly exponential contraction}\label{contraction}
  
  \begin{thm}\label{contraction thm}
Let $f$ and $g$ be two hybrid equivalent quadratic-like maps with modulus at least $\mu$.
Assume that $f$ and $g$ are $n$ times renormalizable. Then
    $$\dist (R^n f, R^n g)\leq C\rho^n,$$
where $\rho\in (0,1)$ is an absolute constant, and $C>0$ depends only on $\mu$. 
\end{thm}

\begin{pf} Let us fix a big $\bar p_e$. By Theorem \ref{puzzle}, if $f$ is a renormalizable map 
with $p_e(f)\geq \bar p_e$,  then the renormalization $Rf$ is $\eps$-close to a quadratic map,
where $\eps=\eps(\mu,\bar p_e)$.
% Let us consider the natural parametrizations of the hybrid
% classes $\HH(f)$ and $\HH(Rf)$ by the space $\EE$ of analytic circle endomorphisms.
% Then $R$ induces a self map $\EE\ra \EE$ (which will be denoted by the same letter). 
It follows that there is an absolute $\delta>0$ such that the Banach ball $\BB_V(f,\delta)$
is mapped into the Banach ball $\BB_U(Rf, \eps)$ (here $V$ and $U$ are the
appropriately chosen domains of $f$ and $Rf$ respectively). By the Schwarz Lemma
in Banach spaces (see the Appendix), this map is uniformly $\rho$-contracting once
 $\eps<{1\over 2}\rho\delta$ 
(which is the case  for sufficiently big $\bar p_e$). 

Let us now show that there is an $N=N(\bar p_e)$ with the following property:
If for $n\geq 2N$ consecutive renormalization iterates,
 the essential period  stays bounded by $\bar p_e$,
then $R^N$ is contracting by 1/2.   Indeed, otherwise we can find a sequence of hybrid equivalent
finite towers $\bold F_n=\{F_{m}\}_{m=-l(n)}^{l(n)}$ and $\bold G_n=\{G_{m}\}_{m=-l(n)}^{l(n)}$ 
 of growing height $l(n)\to \infty$ but such that the $F_0$ and $G_0$ stay a definite
distance apart (where $n$ runs over a certain subsequence of $\N$). Passing to a geometric limit, we
come to a contradiction with the Rigidity Theorem for parabolic towers \ref{Ben}.

Let us now consider the mixed case. Let $\eps>0$.   By Teichm\"uller non-expansion 
(see \S \ref{q-l maps})
and the relation between the Teichm\"uller and Banach metrics (Lemma \ref{Teich metric}),
 a Banach ball in any hybrid class of a sufficiently small radius
$\delta>0$ is mapped by all iterates $R^n$ into a Banach $\eps$-ball.
By the Schwarz lemma, all the iterates $R^m$ have a uniformly bounded norm on the hybrid 
classes.

Now, assume that  we have  $n\leq 2N$ renormalization iterates of essentially bounded type
(i.e., with essentially period bounded by $\bar p_e$) followed by the renormalization of high type
(i.e., with essential  period greater than $\bar p_e$). If the quantifier $\bar p_e$ is selected
to be sufficiently big, then the contraction factor $\rho$ of the last iterate suppresses the
bounded expansion factor of $R^n$. Thus the whole composition is uniformly contracting.

On the other hand, if $n\geq 2N$ then we have a uniformly exponential  contraction of the 
first $n$ iterates. Indeed, 
 every  cascade of $N$ consecutive iterates of essentially  bounded type 
(except perhaps the tail cascade) 
 contracts by $1/2$, while the tail cascade has a priori bounded norm.
Altogether this yields the desired.  
\end{pf}

\comm{subsection{Overflowing and transversal non-linearity}\label{non-linearity}

Given a quadratic-like family $\Bf=\{f_\la\}$, let $\mod(\Bf) =\inf\mod(f_\la)$
and $M_\Bf$ stand for the Mandelbrot set in $\Bf$. 

\begin{lem}\label{overflowing}
Let $\mu>0$, and $\Bf$ be a proper quadratic-like family with winding number 1 and
$\mod(\Bf)\geq\mu$. Then  for any $M\in \MM$, $R_M\Bf$ is a proper quadratic-like family
with winding number 1, and 
$$\mod(R_M\Bf)\geq \nu(\mu, M)\geq\bar\nu(\mu)>0,$$ where
$\nu(\mu, M)\to\infty$ as $p_e(M)\to \infty$.
\end{lem} 

\begin{pf} By \cite{DH:pol-like}, 
The family $R\Bf$ is proper with winding number 1. 
By \cite{puzzle}, Theorems I,II, if $f$ is a renormalizable $\mod(f)\geq \mu>0$, then
$\mod(R\Bf)\geq \nu(\mu)>0$.  By \cite{puzzle}, Theorems III,V,
$\mod(R\Bf)\geq \nu(\mu,M)\to\infty$ as $p_e(M)\to\infty$. 
\end{pf}

\begin{cor}\label{nonlinearity} Under the circumstances of the above lemma,
  $R$ has a bounded nonlinearity on $M_\Bf$ (depending only on $\mu$). 
\end{cor}

\begin{pf}%  Let us consider the set $\XX$ of proper quadratic-like families $\Bf$ with
% winding number 1 and such that
% $\mod(\Bf)\geq\nu/2$, $\Psi(f_\la)=r$. If $r$ is sufficiently small (depending only on  $\mu$)
% then $\XX$ is compact in the open compact topology.
Koebe Theorem.
\end{pf}   
end comm}

\subsection{Cylinders}\label{cylinders sec} 

Let us consider an orbit $\{R^m f\}_{m=-l}^n$. Its  $(l,n)$-{\it itinerary} is the sequence of
the Mandelbrot copies $\{M_m\}_{m=-l}^n$ such that $R^m f\in \TT_{M_m}$.
% The set of points
% with the same itinerary as above is called a {\it cylinder} of rank $(l,n)$.  

\begin{lem}\label{cylinders} 
Let us have two points $f$ and $g$ with the same $(l,n)$-itinerary  and such that
$\mod(R^k f)\geq\mu>0$ and $\mod(R^k g)\geq \mu>0$, $ -l\leq k\leq n$.  
Then $\dist(f,g)<\eps=\eps(\mu, l, n)$, where $\eps\to 0$ as $l,n\to\infty$
($\mu$ being fixed). 
\end{lem}

\begin{pf}
By Theorem \ref{contraction thm}, there exist $\rho\in (0,1)$ and $N$ such that $R^N$ is  
$\rho$-contracting on the foliation $\FF\cap \QL(\mu)$.  

 Let $\chi(f)=P_c$ and $\chi(g)=P_b$. By Corollary \ref{real rigidity}, 
$|b-c|<\delta(n)\to 0 $ as $n\to\infty$, so that $f$ and $g$ lie on the nearby
leaves of the foliation $\FF$. The same is applicable to $f_k\equiv R^k f$ and
 $g_k\equiv R^k g$, $k=-l,\dots, N$.

For any integer $k\in [-l,0]$, let us consider a  map $h_k\in \HH(f_k)$ 
belonging to the vertical fiber via $g_k$, i.e., 
$\pi(h_k)=\pi(g_k)$.  Then $R^N$ $\rho$-contracts the distance between $f_k$ and $h_k$. On the
other hand, by Theorem \ref{parapuzzle},
$R^N h_k$ and $R^N g_k$ belong to the same quadratic-like family of class $\GG_\nu$ 
with $\nu=\nu(\mu)$. 
Since they  have the nearby straightenings,  Theorem \ref{trans qc} implies that
 $\dist(R^N h_k, R^N g_k)<\delta_1(n)\to 0$ as 
$n\to\infty$.

Take an $\eps>0$ and a $\rho'\in (\rho,1)$, and find an $n$ such that
$\delta_1=\delta_1(n)<{(\rho'-\rho)\eps\over \rho+1}$. If $\dist(f_k, g_k)\geq\eps>0$ then
\begin{eqnarray}\label{nearby contraction}
\dist(R^N f_k, R^N g_k)\leq \dist(R^N f_k, R^N h_k)+ \dist (R^N h_k, R^N g_k)\leq\\
  \rho\dist(f_k, h_k)+ \delta_1\leq \rho((\dist(f_k,g_k)+\delta_1)+\delta_1<
                 \rho'\dist(f_k, g_k). \nonumber
\end{eqnarray}
Thus $R^N$ uniformly contracts the distance between the $f_k$ and $g_k$, while it stays greater
than $\eps$. Hence in bounded number of steps  (depending on $\eps$) 
this distance must become less than $\eps$.
\end{pf} 

\subsection{Realization and rigidity of general towers}\label{general towers} 
Let us now prove that any real  combinatorics
$\bar\tau=\{M_k\}_{k=-\infty}^\infty$, $M_k\in \MM_\R$, can be be realized by a unique tower with a
priori bounds. 

\begin{thm}\label{realization} For any real combinatorics $\tau$ there is a unique tower $\bar f$ with
this combinatorics and a priori bounds. Moreover, this tower is real and
$\mod(\bar f)\geq  \nu$ with an absolute $\nu>0$.
\end{thm}

\begin{pf}
By Theorem \ref{a priori bounds}, there is an absolute $\nu>0$ such that  for any infinitely
renormalizable quadratic polynomial  $f=P_c\in \II$,
$R^n f\in \QL(\nu)$, $n=0,1,\dots$.

Let us take a combinatorial sequence $\bar \tau=\{M_k\}$.
 For any $l\geq 0$,  there is a real infinitely renormalizable
 quadratic polynomial $P_l\equiv P_{c_l}$ with combinatorics 
$\tau(P_l)=\{M_{-l},\dots, M_0, \dots\}$. 
 Let $f_{0,l}=R^l P_l$. These are infinitely renormalizable real quadratic-like maps with 
common combinatorics
 $\tau_0=\{M_0,M_1,\dots\}$ and $\mod(f_{0,l})\geq\nu$.  Since the set of such maps is compact,
  we can pass to a quadratic-like  limit $f_0=\lim_{l\to\infty} f_{0,l}$ (along a subsequence)
   with the same properties. 
  
  Let us now do the same thing for every $i\leq 0$. Let $f_{i,l}=R^{l+i} P_l$, and let
$f_i=\lim_{l\to\infty}  f_{i,l}$
  be a limit point. The map $f_i$ is real and has combinatorics $\tau_i=\{M_i,M_{i+1},\dots\}$.
  
   Selecting the above converging subsequences by means of the diagonal process,
  we  construct a sequence of real infinitely renormalizable
  quadratic-like maps $\{f_i\}_{i=-\infty}^\infty$ such that $Rf_i=f_{i+1}$,
  $\chi(f_i)\in M_i$ and $\mod(f_i)\geq\nu$. 
  This sequence represents a real tower $\bar f$ with combinatorics $\bar \tau$ 
and a moduli bound $\nu$.
  
  Thus any real combinatorics $\bar\tau$ is represented by a tower with a priori bounds.
Moreover,  this tower is unique. Indeed, if $\bar f$ and $\bar g$ are two such towers then by
Lemma \ref{cylinders} 
$\dist(f_0, g_0)$ is arbitrary small, so that $f_0=g_0$. 
For the same reason $f_i=g_i$ for any $i$.  
\end{pf}

Let us now state a more general realization and rigidity theorem for one-sided towers.

\begin{thm}\label{one-sided realization} For any real combinatorial past  $\tau=\{M_k\}_{k=-1}^{-\infty}$
and any $c\in [-2, 1/4)$,
 there is a unique tower $\bar f=\{f_k\}_{k=0}^{-\infty}$ with a priori bounds such that 
 $\chi(f_0)=c$ and $f_k\in M_k$ for $k<0$.
 Moreover, this tower is real and
$\mod(\bar f)\geq  \nu(\eps)>0$, provided $c<1/4-\eps$.
\end{thm}

\begin{pf} Theorem \ref{a priori bounds} for real finitely renormalizable
quadratic-like  maps yields the desired statement in
the same way as for two-sided towers.
\end{pf}

\subsection{Renormalization horseshoe}\label{horseshoe sec}
Let us now consider the space $\Sigma$ of all possible real combinatorial types 
$\hat \tau=\{M_k\}_{k=-\infty}^\infty$, 
where the $M_k\in\MM$ are selected arbitrarily  from the
family of real maximal Mandelbrot copies. Supply $\Sigma$ with the weak topology.
 Let $\omega: \Sigma\ra\Sigma$ stand for the left shift on this space.

% Let $\II_\infty$ stand for the set of infinitely renormalizable quadratic-like maps, and
% $$\II_\infty(\mu)=\II\infty\cap \QL(\mu)=\{f\in \II_\infty: \mod(f)\geq\mu>0\}.$$ 

Let us say that an infinitely renormalizable map $f\in \QL$ is 
{\it completely non-escaping} under the renormalization
if the full renormalization orbit $\{R^n f\}_{n=-\infty}^\infty$ is well-defined on $f$
and $\mod(f_n)\geq \mu=\mu(f) >0$, $n\in \Z$. 
% (for a certain sequence of analytic extensions of the renormalization branches $R_M$). 
Note that in this case by Lemma \ref{injectivity}
the backward trajectory $\{R^{-n} f\}_{n=0}^\infty$ 
is uniquely determined by $f$.

Let $\AAA\subset\QL$ stand for the set of completely non-escaping orbits with real 
combinatorics.
We call this set the {\it (full) renormalization horseshoe} for the following reason:

\begin{thm}\label{horseshoe thm}
There exist absolute $\nu>0$ and $\rho\in (0,1)$ with the following properties.
The set $\AAA$ belongs to  $\QL_\R(\nu)$ and 
$R: \AAA\ra\AAA$ is a homeomorphism. There exists  
a homeomorphism  $\eta: \Sigma\ra\AAA$
conjugating $\omega$ and $R|\AAA$. Moreover, for any real infinitely renormalizable map
$f$ there exists a $g\in \AAA$ such that 
\begin{equation}\label{exp conv}
\dist(R^n f, R^n g)\leq C \rho^n,
\end{equation}
 where $C$ depends only on $\mod f$.
\end{thm}

\begin{pf}  
% By Lemma \ref{} any completely invariant map is infinitely renormalizable
Any completely non-escaping point $f$ generates a bi-infinite tower
 $\{R^n f\}_{n=-\infty}^\infty$
with a priori bounds, and vice versa: zero coordinate of such a tower is non-escaping.
By Theorem \ref{realization}, any combinatorics $\tau\in \Sigma$ can be realized
by  unique such a tower.
%Thus the points of $\AAA$ are encoded by bi-infinite towers with a priori bounds. 
Thus we can define a map $\eta:\Sigma\ra \CC(\nu)$ by associating to a combinatorics
$\tau\in\Sigma$ the zero coordinate $f_0$ of the tower $\bar f=\{f_i\}$ representing $\tau$.
 This map is continuous by Lemma \ref{cylinders}. Let $\AAA$ be its image.
 Clearly, $\AAA$ is $R$-invariant and $\eta$ conjugates
the shift $\omega$ and $R|\AAA$. Moreover, by Lemma \ref{injectivity} this map is injective.

Since $\omega$ is a homeomorphism, $R:\AAA\ra\AAA$ is bijective. 
Let us show that it is a homeomorphism. 
Let $J_k=M_k\cap \R$ be the
real traces of the little copies $M_k\in \MM$. Let us denote by $\JJ=\{J_k\}$
 the family of these intervals (formally the same as $\MM$). 

For $J\in \JJ$, let $\TT_{J}=\TT_{M}\cap \QL_\R$ be the corresponding strip of real
quadratic-like maps, and $\AAA_J=\AAA\cap \TT_J$. As the boundary points of
$J$  are exactly once renormalizable, $\AAA\cap\di \TT_J=\emptyset$ for any $J\in \JJ$.
Hence any map $f\subset \AAA_J$ belongs to $\TT_J$ together with 
some neighborhood $\UU$. Since every branch $R_M$ of the renormalization is continuous,
$R$ is continuous at $f$. 

Let us show that $R^{-1}|\AAA$ is also continuous. Let $f\in \AAA$ and $R^{-1}f\in \AAA_J$.
Let $I$ be any other interval of family $\JJ$. Then by Lemma \ref{injectivity}
$R\AAA_I\not\ni f$. Since the strip
$\AAA_I$ is compact, its image $R\AAA_I$ misses some neighborhood of $f$. 

Let us show that these images cannot accumulate on $f$. 
Indeed, otherwise by Theorem \ref{limit points}, $f$ is either parabolic quadratic-like map
or a quadratic polynomial. But in the former  case the map is at most finitely renormalizable,
while in the latter it is not anti-renormalizable 
(see Lemma \ref{little}).

Thus there is a neighborhood $\UU$ of $f$ which misses all the images $R\AAA_I$ with $I\not=J$.
Hence on this neighborhood $(R|\AAA)^{-1}=(R|\AAA_J)^{-1}$. But the latter map is continuous
since $\AAA_J$ is compact. 

Let us now show that $\eta: \Sigma\ra\AAA$ is also a homeomorphism. The only thing to check
is that the inverse map is continuous. Let $f\in \AAA$ be a map with itinerary
$\{J_k\}_{k=-\infty}^\infty$. Let $n\geq 0$. Since $R^k f\in \inter\TT_{J_k}$ for all $k$
and $R|\AAA$ is a homeomorphism, 
all the maps $g\in \AAA$ near $f$ have the same itinerary $(J_{-n},\dots, J_n)$.
But this exactly translates into continuity of $\eta^{-1}$.

Finally,
for any real infinitely renormalizable quadratic-like  map $f$, 
there is a map $g\in \AAA$ with the same combinatorics (by the Realization Theorem
\ref{realization}).  By the Rigidity Theorem \ref{rigidity},
$f$ and $g$ are hybrid equivalent,
 and Theorem \ref{contraction thm}  yields (\ref{exp conv}). 
\end{pf}

\subsection{Periodic points of $R$}\label{periodic points}
It is proven in \cite{universe} that any infinitely renormalizable map $f$ with 
periodic combinatorics is a hyperbolic periodic point for $R$. 
For what follows we will need    
the following weaker statement:

 \begin{lem}\label{non-attracting}
 Periodic points of $R$ are not attracting.
 \end{lem}
 
 \begin{pf}  Let $R^p f_0=f_0$, $P_{c_0}=\chi(f_0)$. 
 If $f_0$ is  attracting than there is  a neighborhood of
 the hybrid class $\HH(f_0)$ attracted to the cycle of $f_0$. 
 By Lemma \ref{escaping}, all the maps in this neighborhood are infinitely renormalizable with
 the same combinatorics as $f_0$.
 In particular, if $c$ is nearby to $c_0$ then $P_c$ is an infinitely renormalizable map
 with the same combinatorics as $P_{c_0}$ contradicting Corollary \ref{real rigidity}.
 \end{pf}
 
 Let $\tD R_\tr$ stand for the tangent action of $\tD R$ in the
one-dimensional quotient bundle $\tT\QL/\tT\FF$ over $\CC$.  
If $R^pf=f$ then the value $\la(f)=\|\tD R_\tr^p(f)\|^{1/p}$ will be called 
{\it mean transversal multiplier} of $f$.
Let 
\begin{equation}\label{mean multiplier}
\bar \la=\inf_p\inf_{f: R^p f=f}|\tD R_\tr^p(f)|^{1/p}
\end{equation}
stand for the \lq{smallest}" mean transversal multiplier of the periodic points of $R$. 
By Lemma \ref{non-attracting}, $\bar \la\geq 1$. If $\bar\la>1$ then we say that
 {\it the periodic points of $R$ are  uniformly hyperbolic} (this term is justified as 
 by Theorem \ref{contraction thm} $R$ is uniformly contracting on the foliation $\FF$). 

\subsection{Invariant cone field and line bundle}\label{cone field}

\begin{lem}\label{derivative bound}
  For any $\mu>0$ and $q\in (0,1)$,  there exist a   $\delta>0$  and a $c>0$   
with the following property. If $f\in \QL_\R$ is $p$ times renormalizable with
$\mod(f)\geq\mu$, then $$\|\tD R_\tr^p(f)\|\geq c(q\bar \la^{\delta})^p,$$
where $\bar \la$ is the smallest mean transversal multiplier (\ref{mean multiplier}). 
\end{lem}

\begin{pf} Let us have two points $f,g\in\CC$ with the modulus at least $\mu$
 lying on the same leaf of the foliation $\FF$. Assume that  $\dist(f,g)\leq\eps$
(which is not necessarily small). 
Then 
\begin{equation}\label{power}
|\tD R_\tr(g)|\geq q\; \|\tD R_\tr(f)\|^{\delta},
\end{equation}
where $\delta>0,\; q\in (0,1)$, 
and $q\to 1$  as $\eps\to 0$. Moreover, there is a $\bar p$ such that if
$p(M)\leq \bar p$ then we can let $\delta=1$, and otherwise we can let $q=1$. 

This follows from the fact that the holonomy from $f$ to $g$ is transversally
qc (Theorem \ref{trans qc}).
% where $K(\eps)\to 1$ as $\eps\to 0$ 
Indeed, (\ref{power}) is obviously true for any particular renormalization $R_M$
with $\delta=1$ and $q=q(M)$, as it is analytic. 
So let us take a Mandelbrot copy $M$ with a big period
$p(M)$. Then by Corollary \ref{shrinking}, $\diam M$ is small. Hence by bounded transversal 
non-linearity (Lemma \ref{nonlinearity}), $\|\tD R_\tr (f)\|$ is big.  

Let us now take some full unfolded quadratic-like families $\SS, \XX\in \GG_\sigma$
via $f$ and $g$ respectively. Let us take a disk $D\subset\SS$ around $f$ of size $\xi>0$
whose image under the renormalization has size of order 1. By transversally bounded
non-linearity, $\|\tD R_\tr(f)\|\asymp \xi^{-1}$. 

Furthermore, since the holonomy $\gamma: \SS\ra\XX$ is qc, it is
H\"older continuous with some exponent $\delta=\delta(\mu)>0$ and an absolute constant.
Hence $\diam (\gamma D)=O(\xi^{\delta}) $, so that 
$$
     \|\tD R_\tr(g)\|\asymp (\diam(\gamma D))^{-1}\geq c \|\tD R_\tr(f)\|^{\delta}.
$$
Finally, since $\|\tD R_\tr(f)\|$ is big, we can kill the constant $c$ by a small decreasing
of the exponent. This yields (\ref{power}). 

Since by \corref{nonlinearity} transversal non-linearity of $R$ is also bounded,
we conclude that the same estimate holds under the assumption that $f$ and $g$ belong
to the same renormalization strip and $\dist(\chi(Rf),\chi(Rg))<\eps$
 (with the constants independent of the strip).

Given an $f$, let us consider the periodic point $g$ of period $p$ which has the same itinerary
$(M_{0},\dots  M_{p-1})$ as $f$. Then by Lemma \ref{cylinders} the orbit of 
$\{R^kg\}_{k=N}^{p-N}$ $\eps$-shadowes the corresponding orbit of $f$, where 
$N=N(\eps)$, and the desired estimate follows from (\ref{power}) by the chain rule.  
%through the renormalization strips. 
%By transversal quasi-conformality we can move $f$ around the leaf (which changes only the
%constant $C$), so that $f$ and $g$ would belong to the same $\pi$-fiber.
%
%Furthermore, let us consider the last moment $l$ when $R^l f$ and $R^l g$ belong to a
% renormalization strip of big period. Until this moment they shadow each other in vertical
%direction as well, so that they have comparable derivatives in the above sense. After that moment
%they still $\eps$-shadow each other except perhaps the last $N=N(\eps)$ iterates.
% But at each of these
% last moments the derivative is bounded, so that this just contributes to the constant.    
\end{pf}

For $f\in \CC_V$, let us define the $\theta$-{\it cone}  as follows;
$$ 
    C_f^\theta=\{u\in \tT_f\QL_V: \theta(u)>\theta\},
$$
where the angle $\theta$ is defined by (\ref{theta})
 (for notational convenience we skip dependence of the cone on the Banach slice). 

\begin{lem}\label{cones}
%Let $\mu>0$ be smaller than the absolute modulus bound $\nu$ from Theorem~\ref{a priori bounds}. 
There exist $\theta>0$, $N$,
% depending on $\mu$,
 and a choice of finitely many Banach slices $\BB_f\equiv \BB_{V(f)}$, $f\in \AAA$, 
such that $R^N C_f^\theta\subset C_{R^N f}^{2\theta}$.
\end{lem}

\begin{pf} 
Let $\nu>0$ be the absolute modulus bound $\nu$ from Theorem~\ref{a priori bounds}.
There is an $\eps>0$ such  that if $f\in \AAA$ and $U$ is a neighborhood of $K(f)$
with $\mod (U\sm K(f))\geq \nu$, then $U$ compactly contains the $(\eps\diam K(f))$-neighborhood 
$\Om_\eps(f)$ of $K(f)$. By compactness, we can find a finite family of fundamental domains
$V(f)\subset \Omega_\eps(f)$, $f\in \AAA$.  By Theorem~\ref{a priori bounds},
there exists an $N$ such that $\mod(R^N f)\geq \nu$. Hence for small enough $\delta>0$, we have:
$R^N \BB_{V(f)}(f,\delta)\subset \BB_{V(R^N f)}$.  

Let us consider the family of tangent cones $C_f^\theta\subset \BB_f$. 

 By Theorem \ref{contraction thm}, $R$ is uniformly exponential contracting in the
$\FF$-direction.
On the other hand, by Lemmas \ref{non-attracting} and \ref{derivative bound},
 $R$ can only slowly contract at the transversal direction:
$\|\tD R_\tr^n(f)\|\geq C q^n$, $f\in \AAA$, with $q$ arbitrary close to 1
(and $C=C(q)$).   

Let us consider the tangent line bundle $\KK=\{K_f\}$ complementing $\tT\FF$ (see (\ref{kernels}).
Recall that $u^h$ and $u^v$ stand for the horizontal and vertical
projections of $u\in \BB_f$.
Since both subbundles $T\FF$ and $\KK$
are continuous, $\|u^v\|\asymp\|u\|_\tr$.

It follows that for  $N$ big enough, there exist $\rho'>\rho>0$ with arbitrary
small ratio $\rho/\rho'$  such that 
\begin{equation}\label{e1}
\|(R^N u)^h\| \leq \rho\|u^h\|,\quad  \|(R^N u)^v\|\geq \rho' \|u^v\|,
\end{equation} 
where the norms are taken  in the corresponding Banach spaces $\BB_f$.

Let now $u\in \di C_f^\theta$ with a small $\theta$, $r\in (\rho, \rho')$. Then 
\begin{equation}\label{e2}
\|(R^N u)^h\|=\|R^N u^h + (R^N u^v)^h\|\leq \rho\|u^h\|+ O(\theta  \|u^h\|)\leq r\|u^h\|,
\end{equation}
provided $\theta$ is sufficiently small. 
By  the second inequality of (\ref{e1}) and (\ref{e2}),  $R^N u\in C_{R^N f}^{2\theta}$.
\end{pf} 

We are ready to construct the unstable-to-be line bundle $\EEE^u$ over $\AAA$. 

\begin{lem}\label{line bundle}
The renormalization operator has a continuous invariant  tangent line field 
$\EE^u=\{E^u_f\subset T_f\QL\}$ over the horseshoe $\AAA$ 
transversal to $\FF$.
% over $\CC(\mu)$.
\end{lem}

\begin{pf} This is a standard construction by going backwards and pushing forward the cones:
Let $C^{\theta,n}_{f}=R^n C_{R^{-n}f}\theta$ and
$E_f=\bigcap_{n\geq 0} C^{\theta,n}_{f}.$ 
Note that since $f\in \AAA$, $R^{-n}f\in \QL(\nu)$ for all $n\geq 0$, 
the cones $C^{\theta,n}_{f}$ are well-defined  and nested by the previous lemma. 

Let us consider the projective cone $\hat C_f^\theta$. i.e., the space of lines in $C_f$. 
It can be realized as the cross-section of $C_f$ by the hyperplane $\{u: u^v=\text{const}\}$. 

Supply the projective cones with the projective distance as follows:
For $\hat u, \hat v\in \hat C_f$, consider the line interval
$I(\hat u,\hat v)=\{w=\hat u+t\hat v\in \hat C_f^\theta\},$ and view it 
as the one-dimensional hyperbolic line  $\Hyp^1$
Then the projective distance between $\hat u$ and  $\hat v$ 
is defined as the hyperbolic distance between
$\hat u$ and $\hat v$ in $I(\hat u, \hat v)$. 

The embedding $C^{2\theta}_{f}\ra C^{\theta}_{f}$ 
uniformly contracts the projective distance on these cones, while the differential
$\tD R^n: C_{R^{-n} f}^\theta\ra C_{f}^{2\theta}$ is at least simply contracting.
Thus $\tD R^n: C_{R^{-n}f}^\theta\ra C_f^\theta$ is uniformly contracting. 

It follows that the projective cones $C^{\theta,n}_{f}$ uniformly exponentially shrink to 
some projective points. These points represent the tangent lines $E^s_f$ transversal to $\FF$.
 This line bundle is clearly invariant. 
It is also continuous: indeed the cone field  $\{C^{\theta,n}_{f}\}$ is continuous for any
given $n$, and well localizes the line field for $n$ big enough.    
\end{pf}

\subsection{Slow shadowing and hyperbolicity}\label{perturbation} 
We will now prove in the similar way as \cite{universe} that  $R|\AAA$ 
is uniformly hyperbolic. The idea is to construct (assuming the contrary)  an 
$\orb (g)$ which slowly shadows some $\orb(f)$ on $\AAA$, which contradicts
the Rigidity Theorem \ref{rigidity}.

 Let $\EE^s$ stand for the tangent bundle to $\FF$ over $\AAA$ 
(the horizontal subbundle) and $\EE^u$ denote as above
the transversal line bundle given by Lemma \ref{line bundle}.
% We will also call the latter \lq{vertical} though it is different from the
% above vertical subbundle.

\begin{thm}\label{hyperbolicity thm}
  The renormalization operator $R: \AAA\ra\AAA$ is uniformly hyperbolic with $\EE^s$ and
$\EE^u$  serving for the stable and unstable subbundles.
\end{thm}

\begin{pf} Due to \lemref{derivative bound}, it is enough to prove uniform hyperbolicity of
the periodic points.  Assume the contrary: $\bar\la=1$.

Let us consider a family of Banach slices $\BB_f\equiv \BB_{V(f)}$, $f\in \AAA$, 
 from Lemma \ref{cones}. By \cite[Lemma 4.12]{universe}, 
the  slice $\FF_f\equiv \FF_{V(f)}$ of the foliation
$\FF$ admits a local extension to a Banach neighborhood $\UU_f\subset\BB_f$ of $f$,
 which will still be denoted as $\FF_f$.  Since $\bar\AAA$ is compact in $\QL$, 
the size of $\UU_f$ is uniform over $f\in \AAA$.
The leaf of $\FF_f$ via $g$ will be denoted as $\LL_f(g)$, $\LL(f)\equiv\LL_f(f)$ .

Let $E^h_f(g)\subset \BB_f$ stand for the tangent plane to the leaf $\LL_f(g)$ of $\FF_f$ via $g$,
and $E^v_f(g)\subset \BB_f$ stand for the complementary line via $g$ parallel to 
$\tT_f \ZZ_f\equiv E_f^v(f)\equiv E_f^v$. Given a tangent vector $u\in \tT_g\BB_f$, its 
projections to $E^h_f(g)$ and $E^v_f(g)$ will be respectively denoted as $u^h$ and $u^v$.
% For $g\in Q_f$ and a vector $v\in \tT_g Q$, 
Let us define the angle $\gamma=\gamma(u)$ by 
$\tg(\gamma)=\|u^v\|/\| u^h \|$.  Let 
$$
   \Lambda_f(g)=\{u\in \tT_g \BB_f:\;  \gamma(u)>\pi/4\}
$$  
 stand for the $\pi/4$-cone at  $g$.
% $\Lambda_f\equiv \Lambda_f^{\pi/4}$.   

Let us consider a topological bidisk $Q_f\subset \BB_f$ centered at $f$ of such a kind.
Take a vertical topological disk $\SS_f\subset E^v_f$ containing $f$, and consider its motion
$\SS_f\ra \SS_f(g)$ under the holonomy along $\FF_f$, where $g$ runs over a neighborhood
$\VV_f\subset \LL(f)$ of $f$.  The disks $\SS_f(g)$ will be called the {\it vertical
cross-sections} of $Q_f$.   
%
% Let us select an analytic Banach local chart near any $f\in \AAA$  modeled on
% $\tT_f\QL=E^s_f\oplus E^u_f$ and depending continuously on $f$. 
% Within this chart, let $P^s: \tT_f\QL\ra E^s$ and $P^u: \tT_f\QL\ra E^u $
% correspond to the stable and unstable projections respectively.
%
%  Let us call a {\it bidisk}
% $Q_f=Q^s_f(\delta)\times Q^u_f$  a product in  the above local chart  of a
% stable Banach $\delta$-disk and an unstable topological disk. 
% An unstable cross-section of such a bidisk is defined as
% $\SS_g=g\times Q^u_f$, where $g\in Q^s_f$.
% and similarly for the horizontal cross-section.
% It is a one-dimensional quadratic-like  family.  
%
Let  $$
   \di^h Q_f=\bigcup_{g\in\VV_f} \di \SS_f(g)\quad \text{and}\quad 
    \di^v Q_f=\bigcup_{g\in \di\VV_f} \SS_f(g)
$$
 stand  respectively for the horizontal and vertical boundaries  of the bidisk $Q_f$.

% Let $$\La^\gamma_f=\{(g,u): g\in D_v, \theta(u)>\pi/4\}$$
%  stand for the $\pi/4$-cone field over $D_f$. 
%
Let $\Upsilon_f $ stand for the family of complex analytic curves  $\Gamma$ in $Q_f$ 
with $\di\Gamma\subset \di^h Q_f$ and 
 with the tangent lines $\tT_g\Gamma$ belonging  to the cones $\La_f(g)$.

We can select the family of bidisks $Q_f$ in such a way that  for some $N$
it satisfies the following properties:
\begin{itemize}
\item {\it Horizontal contraction}: There is a $\rho<1$ such that
if $g\in Q_f$ and $R^N g\in Q_{R^n f}$ then
 for any $u\in \tT_g Q$, $\|(\tD R^N v)^h\|\leq \rho \|v^h\|$;
% $R^N D_f^h\cap D^_{R^Nf}=\emptyset$;
%
\item  {\it Invariance of the cone fields}: If $g\in Q_f$ and $R^N g\in Q_{R^n f}$ then 
  $R^N \Lambda_g\subset \Lambda_{R^n g}$;
 \item {\it Overflowing property for high periods}:
 There exist $\mu>0$ and  $\bar p$ such that if $p(f)\geq \bar p$ 
% and $\Gamma\in \Upsilon_f$,
  then for any $g\in Q_f^s$, $R^N\SS_g$ is a full family of class  $\GG_\mu$ and
  $R^N \SS_g\cap Q_{R^N f}\in \Upsilon_{R^N f}$; 
 \item  {\it Definite vertical size on every strip}:
 For any $M$, there is an $\eps=\eps(M)>0$ such that
  any vertical cross-section $\SS_g$ of $Q_f$ contains a round disk of radius $\eps$ 
  centered at $g$.
\end{itemize}  

To make such a selection, let us first take a $\delta>0$ so small that
 $\Pi \BB_f(f, \delta)\subset \HH_{W(f)}(\Pi(f), \eps)$, $f\in \AAA$, 
where the vertical tube over  $\HH_{W(f)}(\Pi(f), \eps)$ can be equipped with a
holomorphically moving fundamental annulus (see Lemma \ref{equipment}). 

Now, all further selections will be made in such a way that $\diam Q_f<\delta$, 
so that all the bidisks $Q_f$
will belong to the corresponding vertical tubes. 
Let us  select the horizontal section $\VV_f$ of a bidisk $Q_f$ as a Banach
ball in the leaf $\LL_f$ of radius $\delta/4$. 
The choice of the vertical cross-sections
$\SS_f$ depends on the renormalization strip $\TT_M\ni f$ in the following fashion.  

Let $\D_f(g,\eps)\subset E_f^u(g)$ stand for the round disk of radius $\eps$ centered at $g$. 
Let $\MM_f(g)= E_f^u(g)\cap \CC$ stand for the Mandelbrot set in the complex line $E^u_f(g)$
 considered as a quadratic-like family. 
By Corollary \ref{shrinking}, the little Mandelbrot copies $M_i$ in the quadratic
family $\QQ$ shrink.  As  by Theorem \ref{trans qc} the foliation $\FF$ is transversally
quasiconformal,  
 $\MM_f(g)\subset \D_g(\delta/4)$, provided  $p(M)\geq\bar p$ with big enough $\bar p$.
For $p(M)<\bar p$, we let $\SS_f=\D_f(\delta^{K})$, where $K>1$ is  a bound on the qc
dilatation  of the holonomy along $\FF$.  Then
the sections $\SS_f(g)$ obtained from
$\SS_f$ by the holonomy have diameter at most $\delta$ (since $K$-qc maps are H\"older 
continuous with exponent $1/K$).  

For $p(M) \geq \bar p$, let us consider the vertical fiber $\ZZ_f$ via $f$. 
By the above choice, it is a full unfolded equipped quadratic-like family with a definite geometry.
By Theorem \ref{parapuzzle}, the renormalization $R_M$ analytically extends to a parameter
puzzle piece $\tl\Delta_f(g)\subset \ZZ_f(g)$ bounded by a $\kappa$-quasicircle
 and such that  $\mod(\tl\Delta_f(g))\asymp\nu>0$,
with  absolute $\kappa>0$ and $\nu>0$. 
% Moreover, for different $g\in \VV_f$, these puzzle pieces 
% are related by holonomy along the foliation $\FF$. 
Let us shrink these domains a bit, $\Delta_f(g)\Subset \tl\Delta_f(g)$, so that
both $\mod(\tl\Delta\sm \Delta)$ and $\Delta_f\sm M_f$ are still definite,
where $M_f=\ZZ_f\cap \CC$ 
(in order to provide us  with a Koebe space for the renormalization).    

 Let $\SS_f(g)\subset E_f(g)$ be the image of $\Delta_f$ by the holonomy $\ZZ_f\ra E^u_f(g)$. 
It follows that these are  $\rho$-quasidisks with a definite space around the corresponding
Mandelbrot sets: $\mod(S_f(g)\sm \MM_f(g))\asymp\mu>0$, where $\rho$ and $\mu$ are absolute.   
Moreover, the corresponding bidisks $Q_f$ have diameter less than $\delta$ if $\bar p$ is
sufficiently  big. Hence these bidisks belong to the vertical tubes of Lemma \ref{equipment},
so that by Theorem \ref{parapuzzle} (and the Remark afterwards)
the renormalization admits analytic continuation to them.

\comm{ let us consider two boxes
 $Q_f\subset Q_f'$ whose vertical cross-sections represent quadratic-like families
 $\SS\subset\SS'$   intersecting the renormalization strip $\TT_M$ by a little Mandelbrot copy
$M_\SS$ such that the moduli $\mod(\SS\ssm M_\SS)$ and $\mod(\SS'\ssm \SS)$ are of order 1,
and the diameters of $\SS$ and $\SS'$ are comparable with $\diam M_\SS$. 
 Moreover, by Lemma \ref{}, the choice can be made in such  a way that
$R_M$ is well defined on the vertical  cross-sections $\SS'$.

Note that by Corollary \ref{shrinking}, $\diam Q^u_f\leq \eps$ for all $f\in \AAA$,
provided $\bar p$ is sufficiently big. }
%
% The parameters $\delta$, $\bar p$ and $\eps$ are to be specified below.   
Furthermore,
if $\delta>0$ and $\kappa>0$ are sufficiently small then 
there exist $0<\rho<\rho'<1$ with arbitrary small ratio
$\rho/\rho'$, and an $N$ such that for any $f\in \AAA$ and $g\in Q_f$ with 
   $\dist(R^N f, R^N g)<\kappa$,
%
%\begin{equation}\label{first}
%\dist(R^N f, R^N g)\leq \rho\dist (f,g)\; \text{if }\; g\in E^s_f\; \text{and}\;
%                    \dist(f,g)<\delta,
%\end{equation}
%
\begin{equation}\label{second}
\|\tD R^N (v)\|\leq \rho \|v\|, \; v\in E^s_g,
\end{equation}
\begin{equation}\label{third}
\|\tD R^N (v)\|\geq \rho'\|v\|,\; v\in E^u_g.
\end{equation}
These estimates follow from the contraction in the hybrid classes
(Theorem \ref{contraction thm}), almost repelling in the transversal direction
(Lemma \ref{derivative bound}), and uniformly bounded non-linearity on the vertical
cross-sections (according to the choice of the boxes). Indeed, these  
immediately  imply  (\ref{third}), and  also imply (\ref{second})
%the horizontal contraction  throughout the box $Q_f$ 
 by means   
of a simple estimate similar to (\ref{nearby contraction}) and the Schwarz Lemma.

Estimates (\ref{second}) and (\ref{third}) yield the first two desired properties
of the family of boxes: horizontal contraction and invariance of the cone field. 
The last property, definite vertical size on the renormalization strips, 
 is obvious from the definition of the boxes. 
 
Let us check the third property, overflowing. By (\ref{second}) and (\ref{third}), the images
$R^N\SS_g$ of the horizontal cross-sections have  the vertical slope at most 1 in scale $\kappa$, 
and the point $R^N g$ stays distance at least $(1-\rho)\delta$
 from the vertical boundary $\di^v Q_{R^Nf}$.
% But  the vertical size of $Q_{R^n f}$ is at most $\eps$, so that 
Hence the overflowing property is
satisfied, once the vertical size of all boxes  is selected to be smaller than 
some sufficiently  small $\eps>0$.
 
 As in \cite{universe}, let us now consider the fiber action $\bar R$ of the renormalization on 
the space $\AAA\times \QL$ fibered over $\AAA$. 
The above boxes $Q_f$ are naturally embedded into the fibers of this space.
Let $\YY$ stand for the union of the embedded boxes.   

%Given a map $f\in\II_\infty$, let $f_k=R^k f$. 
  For $\tau\in (0,1)$ near 1, let us consider a fiberwise contraction
 $T_\tau:\YY\ra\YY$ (linear  on the fibers in the above local charts). 
Consider the perturbation $L_\tau$ of $\bar R^N$  by
postcomposing $\bar R^N$ with this contraction: $L_\tau=T_\tau\circ \bar R^N$. 
 Assume that  a  periodic point $f=f_\tau$ of period $p$
becomes attracting under this perturbation. Then by \cite[Lemma 2.1]{universe}, 
there is an $l$ and a point
\begin{equation}\label{g1}
g=g_\tau\in\di^u Q_{L^{l}_\tau f}
\end{equation}
  such that
\begin{equation}\label{g}
   g\in  \Omega_0=\{h: L_\tau^{k} h\in Q_{ L^{(k+l)} f},\; k=0,1,\dots, \;
               \text{and}\;   L_\tau^{n} h\to \orb(f)\; \text{as}\; \to \infty\}.
\end{equation}

\comm{
 This is  proved in the same way as in the Small Orbits Theorem of \cite{universe}.
Without loss of generality we can let $l=0$. 
Let $\Omega$ be the connected component of $\Omega_0$ containing $f$.
Assume that it does not touch the horizontal boundary $Q^s_f$. 
Let $\di^h \Omega= \di\Omega\ssm \di^u Q_f$. Because of the uniform horizontal contraction,
\begin{equation}\label{Omega}
  L^p_\tau\di^h\Omega\subset \di^h\Omega.
\end{equation}
 Let us consider
the vertical cross-section $\Gamma_0=\SS_f\cap \Omega$ of $\Omega$, and let  iterate it 
forward: $\Gamma_n=L_\tau^{pn} \Gamma_0$. Because of the second property of
the family of boxes (invariance of the cone fields),
the $\Gamma_n$ are the graphs over their projections to $Q^u_f$ with slope at most 1. By 
(\ref{Omega}), the intrinsic boundary of $\Gamma_n$ belongs to $\di^h\Omega$. 
By the first property of the family of boxes (horizontal contraction),
the $\Gamma_n$ converge to an invariant curve $\Gamma$ with the same properties. 
But then
the map $L^{p}_\tau: \Gamma\ra\Gamma$ is an analytic diffeomorphism of this curve
with an attracting fixed point, which is impossible by the classical one-dimensional
Schwarz Lemma.  
end comm} 

Moreover, the overflowing property implies that
 $L_\tau \di^h Q_f\cap Q_{R^N f}=\emptyset$ if $p(f)\geq \bar p$. 
Hence $f_\tau$ with the shadowing point
$g_\tau$ satisfying (\ref{g1}) and (\ref{g})
  may belong only to finitely many renormalization strips $\TT_M$ with $p(M)\leq\bar p$,
and $\dist(f_\tau, g_\tau)\geq \delta_0>0$. 

Since $R^N$ is transversally non-singular, $\dist(R^N f_\tau, L_\tau g_\tau)\geq \delta_1>0$.
Since the vertical diameter of the $Q_f$ goes to 0 as $p(f)\to\infty$, $ R^N f_\tau$ can 
belong to only finitely many possible strips $\TT_M$.
Repeating this argument for all further iterates, we conclude that for any $k\geq 0$,
 $R^{Nk} f_\tau$ may belong to only finitely many strips $\TT_M$ (depending on $k$). 
It follows that
any limit of the maps $f_\tau$ is infinitely renormalizable.

Since we assume that the periodic points are not uniformly hyperbolic
($\bar\la=1$), for any $\tau\in (0,1)$ there is an attracting periodic point
$f_\tau$ and the corresponding shadowing point $g_\tau.$ 
As by Lemma \ref{compactness lemma} the space $\QL(\mu,\rho)$ is compact, we can pass to limits 
$f=\lim f_{\tau_k}$ and  $g=\lim g_{\tau_k}$ as $\tau_k\to 1$.  As we have just shown,
the first function is infinitely renormalizable.
The second one shadows the first:  $\bar R^{Nk} g\in D_{R^{Nk} f}$, $k=0,1,\dots$,
and hence is non-escaping. 
By Lemma \ref{escaping}, $g$ is infinitely renormalizable as well. 

By the Rigidity Theorem
\ref{rigidity}, $g$ must be hybrid equivalent to $f$. 
But on the other hand, (\ref{g1}) implies that $g$ stays on positive distance from
$\HH(f)$ - contradiction.  
%
% Then we perturb $\bar R$ so that weakly hyperbolic
% periodic points will become attracting. The basin of attraction  of such a point must
% spread all way up to the horizontal boundary of some box round one of the periodic points
%of the cycle. Moreover, this box must have a bounded period, since the boxes with high period
% satisfy the overflowing property.
\end{pf} 

\subsection{Global unstable foliation}\label{unstble foliation}
Let us now show that the unstable foliation of the horseshoe $\AAA$ goes through all
real combinatorial classes except the cusp.

 Let us consider a family of disjoint complex  analytic
curves $\gamma$  in $\QL_V$. Given a map $f\in \QL_V$,
denote the curve passing through $f$ by $\gamma_f$. 
Let us say that such a family of curves is {\it normal}
 if for any $f_0\in \QL_V$  there is a  \lq{horizontal-vertical}" 
local chart $\UU=\UU^h\times \UU^v$  such that for any nearby $f$ the
  curve $\gamma_f\cap \UU$ is a graph of an analytic function  $U^v\ra U^h$. 

\begin{thm}\label{unstable foliation}
Take  any $\eps>0$. 
Then there is a family  $\WW^u=\WW^u_\eps$ of complex analytic leaves
$W^u(f)$, $f\in \AAA$ satisfying the following properties.
\begin{itemize}
\item $R^{-1}| W^u(f)$ is well-defined and $R^{-1} W^u(f)\subset W^u(R^{-1} f)$;  
\item If $g\in W^u(f)$ then 
        $\dist(R^{-n} g, R^{-n} f)\leq C\rho^n$ with absolute $C>0$ and $\rho\in (0,1)$;
\item Every unstable leaf $W^u(f)$ transversally intersects at a single point any 
      hybrid class $\HH_c$ with $c\in [-2,1/4-\eps]$;
\item  The family of the unstable leaves $W^u(f)$ is normal;
\item The renormalization $R$ has uniformly bounded non-linearity on all the leaves;
\item The straightening $\chi: W^u(f)\ra \QQ$ is uniformly quasi-conformal;
\item The real traces $W^u(f)\cap \QL_\R$ of the leaves are pairwise disjoint.  
\end{itemize}
\end{thm}

\begin{pf} Let us consider the family of boxes $Q_f$, $f\in \AAA$,
 constructed in the proof of Theorem
\ref{hyperbolicity thm}. We can now add to the properties of this family listed therein
(horizontal contraction, invariance of the cone field, etc.)  the property
of {\it uniform vertical expansion} (stated similarly to the horizontal contraction).
It follows from the hyperbolicity of the horseshoe and uniformly bounded vertical
distortion.

Let us now take a small number $q\in (0,1)$ and scale all the boxes vertically by this factor.
We obtain a family of boxes $\tl Q_f\subset Q_f$. Let us consider a family $\YY_f$ of 
complex analytic curves $\gamma\subset Q_f$ via $f$ whose tangent lines stay within the
corresponding family of cones. Consider also a similar family $\tl \YY_f$ in $\tl Q_f$ but
with additional assumption that these curves spread over the whole cross-section $\tl Q^u_f$.
If $\gamma\in \tl Y_f$ and $R^{k}\gamma\in \YY_{R^k f}$, $k=0,\dots, l$, 
then $R^k \gamma$ is a curve of
$\YY_{R^kf}$, over a vertical domain $D_k\subset Q^u_{R^k f}$. Moreover, by vertically bounded
distortion, the domains $D_k$ are quasi-disks with bounded shape. By the  vertical expansion,
the size of these quasi-disks exponentially grows with $k$. 

It follows that there exists an $l$ such that $R^k\gamma$ goes outside $Q_{R^k f}$ for some $k\leq l$.
By bounded  shape, at this moment $R^k \gamma$ overflows $\tl Q_{R^k f}$, i.e.,
$R^k\gamma\cap \tl Q_{R^k f}\in \tl Y_{R^k \gamma}.$

Now we can construct the local unstable foliation $\tl W^u_\loc$ in
the usual way by  letting 
$$\tl W^u_\loc(f)=\lim_{k\to \infty} R^k \gamma_{-k},$$
where $\gamma_{-k}$ is an arbitrary curve of $\YY_{R^{-k} f}$, and for $\gamma\in \YY_g$,
 $R\gamma$ is understood as the \lq{cut-off}" iterate $R\gamma\cap Q_{Rg}$.
By the overflowing property just established, these manifolds are spread over the whole
vertical cross-sections $\tl Q_f^u$.  
Note that the horizontal contraction yields the usual property: 
\begin{equation}\label{def2}
  \tl W^u(f)=\{g\in Q_f: R^{-n} g\in Q_{R^{-n}f},\; n=0,1,\dots\}.
\end{equation}

Let now $W^u_\loc(f)=R\tl W^u_\loc(R^{-1}f)$. 
By the overflowing property for high periods stated in the proof of Theorem \ref{hyperbolicity thm},
the vertical sizes of these manifolds are bounded away from 0. It follows that
 this family is normal.

Globalize the $W^u_\loc(f)$ by iterating  them forward.
 Then for some $n$ all the
leaves of  $R^n \WW^u_\loc$ will intersect all the hybrid classes $\HH_c$ with
 $c\in [-2,\; 1/4-\eps]$.  Indeed, by Theorem \ref{one-sided realization}, any combinatorial past
$\tau_-=\{\dots, J_{-2}, J_{-1}, c\}$ with $J_i\in \JJ$ and arbitrary $c\in [-2,\; 1/4-\eps]$ can
be realized by a one-sided tower $\{\dots, g_{-1}, g_0\}$. On the other hand, take
a two-sided tower $\{f_k\}_{k=-\infty}^\infty$ with combinatorics  $\{J_k\}_{k=-\infty}^\infty$
 which has the same combinatorial past as  $\tau_-$. By Lemma \ref{cylinders},  $\dist(f_{-k},
g_{-k})<\eps$ for all 
$k\geq n=n(\eps)$.  Together with (\ref{def2}) and the overflowing property for high
periods this implies that  $ g_{-n}\in W^u_\loc(f_{-n})$, as was claimed. 

Transversality between $\WW^u$ and $\FF$ follows from the corresponding property for the
local unstable  leaves and transversal non-singularity of $R$ \cite[Lemma 5.3]{universe}.
 Uniqueness
of the intersection point  follows from the uniqueness of the one-sided tower with
a given combinatorics (Theorem \ref{one-sided realization}).   

Since $R$ is transversally non-singular and the local unstable foliation is normal, 
the family  $\WW^u=R^n\WW^u_\loc$ of global leaves will also be normal.   
Now bounded non-linearity follows from the Koebe Theorem (compare Lemma \ref{nonlinearity}),
while the bounded dilatation follows from Theorem \ref{trans qc}. Disjointness of the real 
traces of the leaves (or better: disjointness of the intersections $W^u(f)\cap\CC$)
 follows from  Lemma \ref{injectivity}. 
\end{pf} 

\section{\bf Consequences}\label{consequences}

\subsection{Proof of Theorem \ref{measure 0}}
Let us take any infinitely renormalizable parameter value $c\in \II$. 
By Lemma \ref{horseshoe thm}, there is a point $f\in \AAA$ with $\chi(f)=c$. Let 
$I_n(f)=R^{-n} W^u (R^n f)\subset W^u(f)$. 
By Theorem \ref{unstable foliation}, $\diam I_n(f)\to 0$ as $n\to \infty$. 
Moreover, the same theorem implies by means of the standard hyperbolic estimate of 
the distortion that 
\begin{equation}\label{map}
R^n: I^n(f)\ra W^u(R^n f)
\end{equation}
 has a uniformly bounded distortion. 

For $g\in \AAA$, let us consider the interval 
 $L(g)=(\chi| W^u(g))^{-1}(-3/4,\; 1/4-\eps)$ on the unstable manifold $W^u(g)$
consisting of maps with attracting fixed point. Since the straightening
$\chi: W^u(g)\ra (-2,\; 1/4-\eps)$  is uniformly quasi-symmetric 
(by Theorem \ref{unstable foliation}), $\diam L(g)/ \diam W^u(g)\geq \delta>0$ for all $g$.

Let now $S_n(f)=R^{-n} L(R^n f)\subset I_n(f)$. Since the distortion of
(\ref{map}) is  bounded, it follows that $\diam S_n(f)/\diam
I_n(f)\geq\delta_1>$ for all $f$ and $n$. But the maps 
in $S_n(f)$ are only $n$ times renormalizable. 
Hence the set  of infinitely renormalizable maps has definite  gaps  
in arbitrary small scales on $W^u(f)$ near $f$. Using once more that the straightening is
uniformly quasi-symmetric we conclude that the same property holds in the real quadratic
family $(-2,\; 1/4-\eps)$ near $c$.  Thus $c$ is not a density point of $\II$, and
the conclusion follows.
\QED

\subsection{Proof of Theorem \ref{uniformly qs}}
Let $J=J^n_i(\eps)$. As in the proof of Theorem \ref{measure 0} we can find an interval
$I=I^n(f)\subset W^u(f)$, $f\in \AAA$, such that $\chi(I)=J$ (see (\ref{map})). Then 
$$
          \sigma^n|J=\chi\circ R^n\circ\chi^{-1}|I.
$$
As $R^n|I$ has bounded distortion and $\chi$ is uniformly qs, the conclusion follows.
\QED

\subsection{Proof of Theorem \ref{MLC}}
Since by Theorem \ref{unstable foliation}, the  family $\WW^u$ of unstable leaves is normal,
there is a neighborhood $\Omega\subset M_*$ of $[-2,\; 1/4-\eps)$ in the Mandelbrot set
 covered by the straightenings $\chi(W^u(f))$ of all leaves. On the other hand, by Lemma
\ref{shrinking}, the maximal Mandelbrot copies $M\in \MM$ shrink as $p(M)\to \infty$. 
Hence there is a $\bar p$ such that $\chi(W^u(f))\supset M$ for any $f\in \AAA$ and any
$M\in \MM$ with $p(M)\geq \bar p$. 

Take a map $f\in \AAA$ with $\chi(f)=c$. Let 
$$M(f)\equiv M^1(f)\supset M^2(f)\supset\dots\ni f$$
stand for the nest of the Mandelbrot copies in the unstable leaf $W^u(f)$ containing $f$.
We have shown that if $p(R^n f)\geq\bar p$ then  $M(R^n f)\subset W^u(R^n f)$. 
But $M^n(f)=R^{-n}
M(R^n f)$, and the map $R^{-n}$ is contracting on the unstable foliation.  It follows that $\diam
M^n(f)\to 0$, provided there is a  subsequence $n_k\to\infty$ such that $p(R^{n_k} f)\geq\bar p$. 

The Mandelbrot sets $M^{n_k}$ have bounded shape because 
on the unstable foliation the renormalization iterates $R^{-n}$ have bounded non-linearity
and the straightening $\chi$ has bounded dilatation. \QED

\section{Appendix: Complex structures modeled on families of
   Banach spaces}\label{appendix}  
This  Appendix is a brief version of \cite[Appendix 2]{universe} included for the
reader's convenience. We refer to \cite{universe} for the proofs and more details.  
 We assume familiarity with the standard theory of manifolds modeled on Banach spaces (see
e.g., \cite{D-thesis,Lang}).

\subsection{Analytic functions theory in Banach spaces}\label{function theory}
Given a Banach space  $\BB$, let $\BB_r(x)$ stand for the ball of radius $r$ 
centered at $x$ in $\BB$, and $\BB_r\equiv \BB_r(x)$. 

\proclaim Cauchy Inequality. Let $f: (\BB_1,0)\ra (\DD_1,0)$ be a complex analytic
map between two Banach balls. Then $\|Df(0)\|\leq 1$. Moreover, for
$x\in \BB_1$,  $$\|Df(x)\|\leq {1\over 1-\|x\|}.$$

The Cauchy Inequality yields:

\proclaim Schwarz Lemma. Let $r<1/2$ and 
$f: (\BB_1,0)\ra (\DD_r,0)$ be a complex analytic map
between two Banach balls. Then the restriction of $f$ onto the ball $\BB_r$ is contracting:
   $\|f(x)-f(y)\|\leq  q\|x-y\|,$
where $q= r/ (1-r)<1.$

\comm{
The point of the following simple lemma as compared with
Theorem \ref{extension thm} is smoothness of the extension and that the parameter space 
 is allowed to be infinitely dimensional.

\begin{lem}[Local extension]\label{local extension}
Let us have a compact set $Q\subset\C$ and a smooth holomorphic
motion $h_\la$ of a neighborhood $U$ of $Q$ over a Banach domain $(\La,0)$. 
Then there is a smooth holomorphic motion $H_\la$ of the whole complex plane $\C$
over some neighborhood $\La'\subset \La$  of 0 whose  restriction to $Q$
coincides with $h_\la$.
\end{lem}
}

Let us state a couple of facts on the intersection properties between analytic
submanifolds which provide a tool to the transversality results.

Let $\XX$ and $\SS$ be two submanifolds in the Banach space $\BB$ intersecting at point
$x$. Assume that $\codim \XX=\dim \SS=1$. Let us define the {\it intersection multiplicity} $\sigma$
between $\XX$ and $\SS$ at $x$ as follows. Select a local coordinate system
$(w,z)$ near $x$ in such a way that $x=0$ and $\XX=\{z=0\}$. Let us analytically parametrize $\SS$ 
near 0: $z=z(t), w=w(t)$, $z(0)=0, w(0)=0$. 
Then  by definition, $\sigma$ is the multiplicity of the root of $z(t)$ at $t=0$. 

\proclaim Hurwitz Theorem.
Under the above circumstances, let us consider a submanifold
 $\YY$ of codimension 1 obtained by a small perturbation of $\XX$. 
 Then $\SS$ has $\sigma$ intersection points with $\YY$ near $x$ counted with multiplicity.

As usual, a foliation of some analytic  Banach manifold is called analytic  (smooth)
if it can be locally
 straightened by an analytic (smooth) change of variable.   

\proclaim Intersection Lemma. Let $\FF$ be a codimension one complex analytic foliation in a domain
of a Banach space.  Let $\SS$ be a one-dimensional complex analytic submanifold intersecting a leaf
$\LL_0$ of the foliation at a point $x$ with multiplicity $\sigma$. Then $\SS$ has  $\sigma$ simple
intersection points with any nearby leaf.

\begin{cor}\label{transversality criterion}
Under the circumstances of the above lemma, $\SS$ is transversal to $\LL_0$ at $x$ if and
only if it has a single intersection point near $x$ with all nearby leaves. 
\end{cor}

Let $X\subset \C$ be a subset of the complex plane. A {\it holomorphic motion} of $X$ over a
Banach ball $(\BB_1,0)$ is a a family of injections $h_\la: X\ra \C$, $\la\in\BB_1$, with
$h_0=\id$, holomorphically depending on $\la\in \BB_1$ (for any given $z\in X$). The graphs of
the functions $\la\mapsto h_\la(z)$, $z\in X$, form a foliation $\FF$
 (or rather a lamination as it is partially defined)  in $\BB_1\times\C$ with complex
codimension 1 analytic leaves. This is a \lqq{dual}"  viewpoint on holomorphic motions.

Given two complex one-dimensional transversals $\SS$ and $\TT$ to the lamination $\FF$ in
$\BB_1\times\C$,
 we have a partially defined holonomy  $\SS\ra \TT$. We say that this map is locally
quasi-conformal if it admits local quasi-conformal extensions near any $(\la,z)\in\SS$. 

Given two points $\la,\mu\in \BB_1,$  let us define the hyperbolic distance $\rho(\la,\mu)$ in
$\BB_1$ as the hyperbolic distance between $\la$ and $\mu$ in the one-dimensional complex slice 
$\la+t(\mu-\la)$ passing through these points in $\BB_1$.  

\proclaim $\la$-Lemma. Holomorphic motion $h_\la$
 of a set $X$ over a Banach ball $\BB_1$ is transversally locally quasi-conformal.  The local
dilatation $K$  of the holonomy from $(\la,z)\in \SS$ to $(\mu,\zeta)\in \TT$ depends only on the
hyperbolic distance $\rho$ between the  points $\la$ and $\mu$ in $\BB_1$. 
Moreover, $K=1+O(\rho)$ as $\rho\to 0$.

\subsection{Inductive limits}
Let $(\V, \bolshe)$ be a partially ordered set 
% satisfying the following properties:
%
% $V$ has a countable base, i.e., there is a countable subset
% $\W\subset \V$ such that 
%  for any $V_1,V_2\in\V$ there exists a {\it common majorant} 
% $W\in\W$ s.t. $W\bolshe V_1$ and 
% $W\bolshe V_2$). We also assume that $\V$ is dense in itself, i.e., for any $V_1\bolshe V_2$, there
% exist $U\in \V$ such that 
% $V_1\bolshe U\bolshe V_2$. 
%
%
 Recall that such a set is called {\it directed}
if any two elements have a common majorant.
% Let us say that $\V$ is {\it weakly directed}
% if  any two elements $V, \tl V\in \V$ 
% can be linked by a chain $V\equiv V_0, V_1,\dots, V_n\equiv \tl V$ such that every pair
% $V_k, V_{k+1}$, $k=0,1,\dots, n-1$, has a common majorant in $\V$. 
% (This definition is motivated by a notion of a marked germ). In what follows $\V$ is 
% assumed to be weakly directed.
%
We  assume that $\V$ has a countable base, i.e.,  
there is a countable subset $\W\subset \V$ such that any $V\in \V$ has a majorant $W\in \W$.

Let us  have a
family of  Banach spaces $\BB_V$ labeled by the elements ov $\V$. An $\eps$-balls in
$\BB_V$ centered in an $f\in \BB_V$ will be denoted $\BB_V(f,\eps)$.
Elements of the $\BB_V$ will be called \lq{maps}\rq (keep in mind further
applications to quadratic-like maps). 
For every pair $U\bolshe V$, let us have a continuous linear {\it injection} $i_{U,V}:
\BB_V\ra\BB_U$. We assume   the following properties:

\begin{itemize}
\item[C1.] {\it Density}: the image $i_{U,V} \BB_V$ is dense in $\BB_U$;

\item[C2.] {\it Compactness:} the  $i_{U,V}$ is compact, i.e.,
the images $i_{U,V} \BB_V(f,R)$ of balls are pre-compact in $\BB_U$. 
 
\end{itemize}

These properties yield:

\begin{lem}\label{simple stuff}
\begin{itemize}
\item If $U, W\bolshe V$, $f\in \BB_V$, $R>0$,  then the metrics $\rho_U$ and $\rho_W$
induced on  the ball $\BB_V(f,R)$ from $\BB_U$ and $\BB_{W}$ are quasi-isometric,
 i.e. $C^{-1}\rho_M\leq \rho_U\leq C\rho_W$ 
(with the constant depending on all the data specified).

\item Let $U\bolshe V$, and $\phi_i: (\BB_U, \BB_V) \ra (\C,\C)$ be a family 
of linear functionals continuous on the both spaces. Let us consider the common kernels
of these  functionals in the corresponding spaces: 
$\LL_U\subset \BB_U$ and $\LL_V=\LL_U\cap\BB_V$. 
Then $\codim(\LL_U|\BB_U)=\codim(\LL_V|\BB_V)$. 
\end{itemize}
\end{lem}

For any $U\bolshe V$, let us identify any $f\in \BB_V$ with its image $i_{U,V} f\in \BB_U$
and span the equivalence relation generated by these identifications. Thus $f\in \BB_U$
and $g\in \BB_V$ are equivalent if there is a common majorant $W\bolshe (U,V)$ such that
$i_{W,V} f= i_{W,U} g$ (then by injectivity this holds for any common majorant). 
% there is a linking sequence
% $U=U_0, U_1,\dots, U_n=V$ with respective common majorants $W_k$ and a sequence of maps
% $f_k\in \BB_{U_k}$, $f_0=f$, $f_n=g$, such that 
% $i_{W_k,U_k} f_k= i_{W_k,U_{k+1}} f_{k+1}$.    
The equivalence classes will be called {\it germs}.
The space of germs is called the {\it inductive limit} 
of the Banach spaces $\BB_V$ and is denoted by $\BB=\underset{\ra}\lim\BB_V$.  
% If $\V$ is directed then $\BB$ is called the {\it inductive limit}. 

Every space $\BB_V$ is naturally injected into the space of germs, and will be considered
as a subset of the latter.
% (with a different intrinsic topology, though).
Given a subset $\XX\subset\BB$, the intersection  $\XX_V\equiv \XX\cap \BB_V$ will be called 
a {\it (Banach) slice} of $\XX$.  

Let us supply it  with the {\it inductive limit topology}. In this topology,
a set $\XX\subset\BB$ is claimed to be open if all its Banach slices $\XX_V$ are open.
The axioms of topology are obviously satisfied. and the linear operations are
obviously continuous (note that the product topology on $\BB\times\BB$ coincides
with the natural inductive limit topology). Thus $\BB$ is a topological vector space.  
Since points are obviously closed in this topology, $\BB$ is Hausdorff.
The following lemma summarizes some useful general properties of inductive limits. 

\begin{lem}\label{properties}
\begin{itemize}
\item[(i)] In the inductive limit topology,
$f_n\to f$ if and only if all the maps $f_n$ and $f$ belong to the same Banach slice 
$\BB_V$ and $f_n\to f$ in the intrinsic topology of $\BB_V$. Any cluster point $f$ of a set
$K\subset\BB$ is a limit of a sequence
$\{f_n\}\subset K$.

\item[(ii)] A set $\XX\subset \BB$ is open if and only if it is sequentially open.

\item[(iii)] If $X$ is a metric space and $\phi: (X,a)\ra (\QL,g)$ is a continuous map  then
there is  neighborhood $D\ni a$ and an element $V\in\V$ such that $\phi D\subset \BB_V$.

\item[(iv)] A set $\KK\subset\BB$ is compact if and only if it is sequentially compact. Such a
set sits in some Banach space $\BB_V$ and bears an induced \lq{Montel metric}"  which is
well-defined up to a quasi-isometry.

\item[(v)] A map $\phi: \BB\ra T$ to a topological space $T$ is continuous  iff every
restriction $\phi|\BB_V$ is continuous. 
%
%\item A map $\phi: \BB\ra \DD$ between two inductive limits (over $\V$ and $\W$ respectively
% is continuous if and only if  for any $V\in\V$ there is a $W\in \W$ and a neighborhood
% $D\subset\BB_V$ such that $\phi \XX\subset\BB_W$ and the restriction $\phi|D$ is continuous.
% The map $ \phi$  is continuous if and only if it is sequentially continuous.
 
\end{itemize}
\end{lem}

{\it Remarks.} 1.
 Any continuous curve $\gamma: (\R,0)\rightarrow (\BB,g)$  locally sits in some space
    $\BB_V$; \smallskip

2. Given a continuous transformation $R: (\BB,f)\ra (\TT,g)$  between
two  spaces of  germs over $\V$ and $\U$ respectively, 
for any $V\in \V$ there exist an $\eps>0$ and an element $U\in\U$ such that
$R  (\BB_{V}(f,\eps))\subset \BB_{U}$. \smallskip

3. The  third statement of the above lemma shows that the space  $\BB$ is not  Fresche, i.e.,
it is not metrizable, and thus does not have a local countable base of neighborhoods.
 However, as we see, the sequential description of basic  topological properties 
(cluster points, compactness, continuity etc.) is adequate in this space. \smallskip  

4. Note that the Banach slices $\BB_V$ are dense in the space of germs $\BB$,  so that their
intrinsic  topology is not induced from $\BB$. However, by compactness, the intrinsic topology on
the Banach balls $\BB_V(f,R)$ is induced from $\BB$.  \smallskip

5. If $\V$ is directed then the space of germ is clearly a linear vector space. \smallskip

Let us define 
a {\it sublimit} of the directed family $\BB_V,\, V\in \V$ as the inductive limit of Banach
spaces
$\BB_V$ corresponding to a directed subset $\U\subset \V$ (which is not necessarily exhausting).

All linear operators  $A: \BB\ra\TT$ between spaces of germs are assumed to be continuous. 
Let us supply this space with the following convergence:
A sequence of linear operators $A_n: \BB\ra\TT$ converges to an operator $A$ 
if for any $V\in\V$ there is a $U\in\U$ such that for all sufficiently big
$n$, $A_n(\BB_{V})\subset \TT_{U}$ and $A_n|\BB_{V}\to A|\BB_{V}$ in the 
uniform operator topology.
%Continuous dependence of an operator on some parameters
%will be understood in the sequential sense.

\subsection{Analytic maps}
We will give most of the definitions in the complex analytic category
automatically accepting the corresponding smooth notions. So the Banach spaces under
consideration are assumed to be over $\C$.

Let us consider an inductive limit $\BB$ over  $\V$.   
By definition,  a function $\phi:\BB\ra\C$ is complex analytic if
all the restrictions $\phi|\TT_V$ are complex analytic in the Banach sense.
 
Let us have a continuous map $R: \VV\ra\BB'$, where $\VV$ is an open subset of $\TT$ and $\BB'$
is an inductive limit space over $\V'$.
 It is called {\it differential} at a point
$f\in\BB$ if there is a real linear operator $A\equiv DR(f): \BB\ra\BB' $ such that
$$R(f+h)-R(x)=Ah+\omega(h),$$
where $\|\omega(h)\|=o(\|h\|)$ in the induced metric  (note that this
makes sense as by \lemref{properties}, $\omega(h)$ locally sits in the ball of a
Banach slice $\BB_V'$ which bears the well-defined, up to quasi-isometry, induced metric).

As usual, a map $R:\VV\ra\BB'$ is called {\it smooth} if it is differentiable at every
point $f\in\VV$ and the differential $DR(f)$ depends continuously on $f$
(which amounts to differentiability of all Banach restrictions).
A map  $R:\VV\ra\BB'$ is called {\it analytic} if it is smooth, and  the differentials
$DR(f)$ are  linear over $\C$. 

\subsection{Complex structures} 
% Let us now have a partially ordered set $\V$  (not necessarily directed)
% with a countable base $\W=\{W_i\}$.

 Let us have a family of Banach spaces $\BB_V$ labeled by elements $V$ of some set $\V$,
and open sets $\UU_V\subset \BB_V$.
% together with inclusions $i_{U,V}$ (for $U\succ V$)
% satisfying  properties C1 and C2. 
%
 Let us have a  set $\QL$ and  
a family of injections $j_V: \UU_V\ra \QL$. 
% compatible with the injections $i_{U,V}$.
 The images  $\SS_V\equiv j_V \UU_V$  are called Banach slices in $\QL$.
The images $j_V \VV_V\subset \SS_V$ of open sets $\VV_V\subset\UU_V$ are called
{\it Banach neighborhoods}. 
We assume the following properties (compare with C1 and C2):
\begin{itemize}
\item[P1:] {\it countable base and compactness.}
 There is a countable family of slices $\SS_i$ with the following
  property: 
Any $f\in \QL$ has a Banach neighborhood $\VV_V$ compactly contained in some $\SS_i$. 

\item[P2:] {\it analyticity.} If some Banach neighborhood $\VV_V\subset \SS_V$ 
  is also contained in another slice $ \SS_U$, then
    the transit map $j_{U,V}= j_U^{-1}\circ j_V :\VV_V\ra \SS_U$ is analytic.

\item[P3:] {\it density.} 
The differential $D j_{U,V}(f) $ of the above transit map has a dense
  image in $\BB_U$.  
\end{itemize}

We endow $\QL$ with the finest topology which makes all the injections $j_V$ continuous
by declaring a set $\VV\subset\QL$ open
if and only if  all its Banach slices $j_V^{-1} \VV$ are open. 
Lemma \ref{properties} should be minor modified in this more
general situation:

\begin{lem}\label{properties1}
In $\QL$,
$f_n\to f$ if and only if the sequence $\{f_n\}$ sits in a finite union of the Banach
slices, and the corresponding subsequences converge to $f$ in the Banach metric.
All other statements of Lemma \ref{properties} are valid in $\QL$ as well with the
modification that  a single Banach slice
 in (iii) and (iv) should be replaced with a finite union of Banach slices.
\end{lem}

We say that a topological space $\QL$ as above is endowed with
{\it  complex analytic structure}
 modeled on the family of Banach spaces. A subset $\QL^\#$ will be called
a {\it slice} of $\QL$ if it is a union of some family  of Banach
neighborhoods $j_V\VV_V$. It naturally inherits from $\QL$ complex analytic structure.   

%\smallskip{\it Remark.} It is certainly not the most general definition one can imagine
% but it is sufficient for our goals. \smallskip
% 
% if the following property is satisfied: If there is a neighborhood 
% $\UU\subset \BB_U$ such that $j_U \UU\subset j_V\BB_V$  then  the transit map 
% $j_V^{-1}\circ  j_U: \UU\ra \BB_V$ is analytic (in particular, this is  the case when 
% $V\succ U$).  

% Any continuous curve $((-\eps, \eps),0)\ra (\QL,f)$ locally sits in some Banach slice. 
% Since the transit maps between the Banach slices  are analytic, we have a well-defined
% notion of a smooth/analytic curve.
By definition, a smooth curve  in $\QL$ locally sits in some Banach slice and is smooth
there.  Since the transit maps between the Banach slices  are analytic, this notion is
well-defined.  
Moreover, the tangency
of two  smooth curves $\gamma_1$ and $\gamma_2$ through $f$ is well-defined via the local
Banach charts as well. Thus we can define a tangent vector to $\QL$ at $f$  as a class of
tangent curves. This is generally not a linear space but rather a union of 
Banach spaces $\tT_f\SS_V\approx \BB_V$
(the space of smooth curves via $f$ lying in $\SS_V$).

Let us call a point $f\in \QL$  {\it regular} if any two Banach neighborhoods
 $\UU\subset \SS_U$ and $\VV\subset \SS_V$ around $f$ are contained in a common slice
$\SS_W$.  At such a point the tangent space $\tT_f\QL$ is a linear space
identified with the inductive limit of the Banach spaces,
             $$\tT_f\QL=\underset{\ra}\lim_{U: f\in\SS_U} \tT_f\SS_U.$$

A map $R: \QL^1\ra\QL^2$   is called {\it analytic} 
if it locally transfers any Banach slice $\SS_U$  to some slice $\SS_V$, and
its Banach restriction     $ j_V^{-1} R\circ j_U$  is analytic. An analytic map has a
well-defined differential $\tD R(f): \tT_f\QL^1\ra \tT_{Rf}\QL^2$ 
continuously depending on $f$ whose Banach restrictions are linear.

An analytic map is called {\it immersion} if it has an injective differential.
The image $\XX$ of an injective immersion $i: \MM\ra\QL$ is called an {\it immersed submanifold}.
It is called an {\it (embedded) submanifold} if additionally $i$ is a homeomorphism onto
$\XX$ supplied with the induced topology. 
For example,  if there is an analytic projection
 $\pi: \QL\ra\MM$ such that $\pi\circ i=\id$ then  $\XX$ is a  submanifold in $\MM$.
(Note that in this case $\pi$ is a submersion at every point of $\XX$.)
By definition, the dimension of $\XX$ is equal to the dimension of $\MM$.  

If $i: (\MM, m)\ra (\XX,f)\subset (\QL,f)$ is an immersion,
then the tangent space $\tT_f\XX$  is defined as the image 
of the differential $Di(m)$.   If the points $m$ and $f$ are regular then $\tT_f\XX$
is a linear subspace in $\tT_f \QL$, so that we have a well-defined notion
of codimension of $\XX$ at $f$. Moreover, if $\MM$ is a
Banach manifold  (in particular, a finite dimensional manifold) then 
$\XX$ locally sits in a Banach slice of $\QL$.

\begin{lem}\label{codim} Let $\NN\subset\MM$ be a connected submanifold in $\MM$.
  Then  $\codim_g\NN$ is constant.
\end{lem}

As usual,
two submanifolds $\XX$ and  $\YY$  in $\MM$
are called {\it transversal} at a point $g\in \XX\cap\YY$ if
$\tT_g\XX\oplus \tT_g\YY=\tT_g\MM$.

\end{document}